\documentclass[10pt, a4paper]{article}
\usepackage[utf8]{inputenc}
\usepackage[english]{babel}
\usepackage{graphicx} 
\usepackage{xcolor}
\usepackage{amsfonts}
\usepackage{mathtools}
\usepackage{hyperref}
\usepackage{amsmath}
\usepackage{bm}
\usepackage{amsthm}
\usepackage{csquotes}
\usepackage{amssymb}
\usepackage{comment}
\usepackage{esint}

\usepackage[paper=a4paper,left=25mm,right=25mm,top=25mm,bottom=25mm]{geometry} 
\newtheorem{thm}{Theorem}[section]
\newtheorem{lemma}[thm]{Lemma}
\newtheorem{proposition}[thm]{Proposition}
\newtheorem{corollary}[thm]{Corollary}
\newtheorem{remark}[thm]{Remark}

\theoremstyle{definition}
\newtheorem{definition}[thm]{Definition}
\usepackage[maxbibnames=9, giveninits=true, doi=false, url=false, isbn=false]{biblatex}
\usepackage{doi}
\usepackage{mathrsfs}
\numberwithin{equation}{section}

\newcommand{\MK}[1]{{\color{orange} [MK: {#1}]}}
\newcommand{\AC}[1]{{\color{violet} [AC: {#1}]}}

\newcommand{\eps}{\varepsilon}
\newcommand{\T}{T}

\newcommand{\R}{\mathbb{R}}
\newcommand{\Z}{\mathbb{Z}}
\newcommand{\N}{\mathbb{N}}
\newcommand{\X}{\mathbb{X}}
\newcommand{\supp}{\textup{supp}} 
\newcommand{\diam}{\textup{diam}} 
\newcommand{\abs}[1]{\left|#1\right|}
\newcommand{\norm}[1]{\left\lVert#1\right\rVert}

\newcommand{\Wper}{W_{\text{per}}}
\newcommand{\Hdots}{\dot{H}^{-s}}
\def\les{\lesssim}
\def\ges{\gtrsim}
\def\lt{\left}
\def\rt{\right}
\def\A{\mathcal{A}}
\def\Tor{\mathbb{T}}
\def\M{\mathcal{M}}
\def\LM#1{\hbox{\vrule width.2pt \vbox to#1pt{\vfill \hrule width#1pt
			height.2pt}}}
\def\LL{{\mathchoice {\>\LM7\>}{\>\LM7\>}{\,\LM5\,}{\,\LM{3.35}\,}}}
\def\restr{{\LL}}

\title{New dimensional bounds for a branched transport problem}
\author{Alessandro Cosenza, Michael Goldman, Melanie Koser}

\begin{document}

\maketitle
\begin{abstract}
    We consider a branched transport problem with weakly imposed boundary conditions. This problem arises as a reduced model for pattern formation in  type-I superconductors. For this model, it is conjectured that the dimension of the boundary measure is non-integer. We prove this conjecture in a simplified 2D setting, under the (strong) assumption of Ahlfors regularity of the irrigated measure. This work is the first rigorous proof of singular  behavior for irrigated measures resulting from minimality.
\end{abstract}

\section{Introduction} 
 Ginzburg-Landau theory \cite{landau} predicts two types of superconductors, usually named type-I and type-II. While in type-II superconductors vortices appear \cite{bbh,sandierserfaty} (that is, singularities of codimension $2$), in type-I superconductors a positive surface energy penalizes the interface between normal and superconducting regions. Such superconductors, when subject to an applied external magnetic field, exhibit complex patterns \cite{pro1,pro3,pro2}. This behavior suggests that branching patterns should appear inside the sample. Results in this direction in the form of scaling laws were first proved for some simplified models in \cite{ChCoKoOt,Choksi2004} and were more recently extended to the full Ginzburg-Landau model in \cite{CoOtSe2016}. Besides the anticipated uniform branching regime where, as its name suggests, the magnetization is (approximately) uniform at the boundary of the sample, a surprising finding of \cite{ChCoKoOt,Choksi2004,CoOtSe2016} is the presence of a regime where it is energetically favourable for the magnetic field to be far from uniform at the boundary. \\ 
 The aim of this paper is to better  understand this non-uniform branching regime and more specifically the magnetization pattern at the boundary of the sample. Due to the various lengthscales involved, the full Ginzburg-Landau model is too complicated to analyze, see \cite{CGOS}.  We thus make a few simplifications.   First, we consider the asymptotic regime of very small external magnetic field in the cross-over between uniform and non-uniform branching regimes. In this regime, formal computations suggest that the normal region collides into a one-dimensional structure. The magnetization, which lives inside the normal region by the Meissner condition, can then be represented as a measure. The reduced model to study should be a branched transport problem inside the sample complemented by a penalization of the deviation from uniformity at the boundary. Because of the magnetic field energy outside the sample, this penalization comes through a Sobolev norm of negative exponent (thus in a weak norm). Let us point out that deep inside the uniform branching regime, a similar model has been rigorously derived in  \cite{CGOS}. In that case, the internal energy remains the same but the boundary condition is rather imposed as a hard constraint.  The rigorous derivation of the model considered in this paper is work in progress. A second simplification we will make in this paper is to consider a 2D model rather than the original  3D model. Actually, since vertical and horizontal variables play a different role, we rather think of our 2D setting as being a one-dimensional but time dependent branched transport problem.  As in \cite{minimizers2d}, where a similar 2D model for the uniform branching regime was studied, the simple geometry of $\R$ as opposed to $\R^2$ turns out to be very useful. In particular, the no-loop property of minimizers, see Lemma \ref{prop:noloop}, implies in dimension one that the branched transport plan is monotone and thus also optimal for the standard optimal transport problem, see Lemma \ref{lem:monotonecoupling}.
\\

 In order to present our main results, let us briefly introduce  the model (see Section \ref{subsec:energydef} for a precise definition). Motivated by the discussion above, we give the definition in arbitrary space dimension.
Let $s\in (0,1)$, $d\in \N$ and $\T>0$. Let $\Tor^d_T=\Tor^d\times(-\T,\T)$ with $\Tor^d=(\R/\Z)^d$, be a sample with thickness $2\T$ and periodic boundary conditions. We represent the normal phase, or equivalently the magnetic field as a measure $\mu=\mu_t\otimes dt$ where, for a.e. $t\in(-\T,\T)$, $\mu_t=\sum_i\varphi_i\delta_{X_i}$, with $X_i\in \Tor^d$, $\varphi_i>0$ and $\sum_i\varphi_i=1$.
We consider for $-\T\leq a<b\leq \T$ the internal energy
\begin{equation}
\label{e:energyintr}
I(\mu, (a,b)) = \int_a^b \sum_{i}\left[\varphi_i^{\frac{d-1}{d}}(t)+ \varphi_i(t) \vert \dot{X}_i(t)\vert^2 \right]\, dt,
\end{equation}
where $\dot{X_i}(t)$ is the derivative of the map $t\mapsto X_i(t)\in \Tor^d$. The first term in \eqref{e:energyintr} can be seen as a perimeter term that favors mass travelling together, while the second term penalizes the kinetic energy (it is a transport term). Hence, this functional is a (non-standard) branched transport energy, see for instance \cite{pegonequivalence,pegonfractal2018,landscape,pegon2023optimal,branchedbook,branchedwellposedness,branchedxia,brancoliniwirth}.  \\
Let us now define the full energy 
\begin{equation}
    E_{s, \T}(\mu, (-T,\T)) = I(\mu, (-T,\T)) + \Vert \mu_{\pm\T} - 1\Vert_{\dot{H}^{-s}}^2,\label{def:energy_s_intr}
\end{equation} 
where $\Vert \mu_{\pm\T} - 1\Vert_{\dot{H}^{-s}}^2=\Vert \mu_{\T} - 1\Vert_{\dot{H}^{-s}}^2+\Vert \mu_{-\T} - 1\Vert_{\dot{H}^{-s}}^2$. The physical case  corresponds to $d=2$ and $s=1/2$, see in particular \cite{CGOS,dephilippis2023energy}, while we mostly focus here on the case $d=1$ but with arbitrary $s\in(0,1)$.\\

Our aim is to study the dimension of the irrigated measures $\mu_{\pm T}$. By symmetry we will focus on $\mu_T$. Based on the construction of low energy states, it was conjectured in \cite{private} that for $d=2$ and $s=1/2$ we should have $\dim \mu_{T}=8/5$, see also \cite{dephilippis2023energy} for a heuristic explanation. In \cite{dephilippis2023energy} it was proven that this conjecture may be reduced to the proof of local scaling laws matching the global scaling laws. Let us recall that by global scaling law we mean the existence of $\beta>0$ such that for a minimizer $\mu$,
\begin{equation}\label{defglobscale}
 E_{s,T}(\mu, (-T,T))\sim T^\beta.
\end{equation}
Similarly, by a local (upper) scaling law we mean a statement of the form
\begin{equation}\label{deflocscale}
 I(\mu,(T-\eps,T))\les_\mu \eps^{\beta} \qquad \textrm{for }  \eps\in(0,T).
\end{equation}
By \cite{ChCoKoOt,Choksi2004,CoOtSe2016}, \eqref{defglobscale} holds with $\beta=3/7$ so that by \cite{dephilippis2023energy}, in order to prove the conjecture from \cite{private} it is enough to prove that \eqref{deflocscale} holds for $\beta=3/7$, see also Remark \ref{rem:critical}. The main insight of this paper is that it is probably too optimistic to hope for a direct proof of a sharp upper bound of the  form \eqref{deflocscale} and that \cite{dephilippis2023energy} contained only half of the puzzle. Calling $\overline{\alpha}$ the conjectured dimension, we now believe that (see Section \ref{sec:dim} for the definition of the dimension of a measure):
\begin{itemize}
 \item   $\dim \mu_T \le \overline{\alpha}$ should result from  the combination of an interpolation argument reminiscent of the proof of the  lower bound in \eqref{defglobscale}, see \cite[Theorem 3.7]{dephilippis2023energy}, together with a first variation argument, see Theorem \ref{thm:lower_local_bound}.
 \item $\dim \mu_T\ge \overline{\alpha}$ should instead result from a combination of a construction, reminiscent of the argument for the upper bound in \eqref{defglobscale}, see Theorem \ref{thm:local_energy_lower}, together with another first variation argument, see \cite[Corollary 3.6]{dephilippis2023energy} and Lemma  \ref{thm:regularity_boundary_measure}.
\end{itemize}
Unfortunately, and contrary to \cite{dephilippis2023energy}, our argument is so far only effective for $d=1$. In this case, based on  considerations analog to  the case $d=2$, see Remark \ref{rem:critical}, the conjectured dimension is
\begin{equation}
\label{e:alphabar}
    \overline{\alpha}(s)=\min\lt(\frac{4}{3}(1-s), 1\rt).
\end{equation}
Notice that this is coherent with \cite[Section 4]{dephilippis2023energy}, where it is proven that in the case of a Wasserstein penalization (which morally behaves like a $H^{-1}$ penalization), $\mu_T$ is  purely atomic.\\
We further assume that $\mu_T$ is Ahlfors regular. For $\alpha \in [0,d]$, we say that a  measure $\sigma$ is upper $\alpha$-Ahlfors regular if for some $M>0$,
\[ \sigma(B_r(x)) \leq Mr^\alpha \quad \textup{ for } x\in \supp \,\sigma \, \textup{ and }\, r\in (0,1).\]
Similarly we say that it is  lower $\alpha$-Ahlfors regular if
\[ \sigma(B_r(x)) \geq Mr^\alpha \quad \textup{ for } x\in \supp \, \sigma \, \textup{ and }\,  r\in (0,1).\]
If both conditions hold (possibly for different $M,m>0$ respectively), we say that $\sigma$ is $\alpha-$Ahlfors regular.


Our main result is the following.
\begin{thm}
\label{thm:Malphareg}
    Let $d=1$, $s\in (0,1)$ and $\mu$ be a minimizer of $E_{s,T}$. If $s\le 1/4$ then (see Section \ref{sec:dim} for precise definitions) $\underline{\textup{dim}}\,\mu_T=1$. If instead $s>\frac{1}{4}$ and   $\mu_T$  is upper $\alpha$-Ahlfors regular, then
    \begin{equation}
    \label{e:alphaless}
        \alpha\leq \Bar{\alpha}
    \end{equation}
    whereas if $\mu_T$  is lower $\alpha$-Ahlfors regular, then 
        \begin{equation}
        \label{e:alphamore}
        \alpha\geq \Bar{\alpha}.
    \end{equation}   
\end{thm}
\begin{remark}
    Theorem \ref{thm:Malphareg} implies in particular that if $\mu_T$ is $\alpha$-Ahlfors regular for $\alpha\in (0,1)$, then 
    \begin{equation}
       \textup{dim}\,\mu_T=\Bar{\alpha}.
    \end{equation}
    Another important implication of Theorem \ref{thm:Malphareg} is that $\mu_T$ cannot be smooth i.e. we cannot have $\mu_T=f dx$ for some $f\in L^\infty$.
\end{remark}

\begin{remark}
    The fact that for $s\le 1/4$, $\underline{\textup{dim}}\,\mu_T=1$ is a combination of Lemma \ref{thm:regularity_boundary_measure} together with the general upper bound Proposition \ref{prop:scaling13} yielding that \eqref{deflocscale} holds with $\beta=1/3$. It thus holds for every $d\ge 1$.
    \end{remark}

    \begin{remark}
    For a positive measure $\sigma$, let us introduce the quantities
\begin{align*}
         &\overline{\textup{dim}}_M\,\sigma = \inf\{\alpha\in [0,d]:\sigma \textup{ is lower }\alpha\textup{-Ahlfors regular}\}, \\
         &\textup{dim}_F\,\sigma = \sup\{\alpha\in [0,d]:\sigma \textup{ is upper } \alpha\textup{-Ahlfors regular}\}.
    \end{align*}
    In the language of \cite{dimension_measure}, these are the upper Minkowski dimension and the Frostman dimension of $\sigma$.
    Theorem \ref{thm:Malphareg} may be rephrased as saying that if  $s>\frac{1}{4}$ and $\mu$ is a minimizer of $E_{s,T}$, then we have (unconditionally)
        \begin{equation}
            \textup{dim}_{F}\, \mu_T\leq\Bar{\alpha}\leq \overline{\textup{dim}}_{M}\, \mu_T.
        \end{equation}
\end{remark}
Let us observe that \eqref{e:alphabar}  as a conjecture is similar in spirit to some questions studied in \cite{pegonfractal2018, pegon2023optimal}. In \cite{pegonfractal2018} a "unit ball" for a branched transport problem is studied and its boundary is conjectured to have a fractal behavior. Similarly, in \cite{pegon2023optimal} the concept of Voronoi cells for branched transport is introduced. Again, interfaces between cells are conjectured to have fractal behavior. It would be interesting to see if the ideas developed here could prove useful also in that context. \\

Let us now comment on the proof of Theorem \ref{thm:Malphareg}. As explained above, the proof of  \eqref{e:alphaless} is a combination of an interpolation argument related to the proof of the  lower bound in the global scaling law \eqref{defglobscale} together with a first variation argument. The interpolation argument can be found in \cite{dephilippis2023energy} and is reported in the upper bound in \eqref{e:dimbound}. The first variation argument is contained in the proof of Theorem \ref{thm:lower_local_bound} below. To prove \eqref{e:alphamore} we combine instead the first variation argument from \cite{dephilippis2023energy} (which we extend here to arbitrary $d$ and $s$ in Lemma \ref{thm:regularity_boundary_measure}) together with a construction which is reminiscent of the upper bound construction for the global scaling law \eqref{defglobscale}. This construction is contained in Theorem \ref{thm:local_energy_lower}.  For both Theorem \ref{thm:lower_local_bound} and Theorem \ref{thm:local_energy_lower}  we need to assume $d=1$.  However since large portions of the proofs are valid for arbitrary dimensions, we decided when possible to keep the presentation general.  Let us now comment on the two bounds separately. The first result that we prove, in order to prove \eqref{e:alphaless}, is the following local scaling law.
\begin{thm}\label{thm:lower_local_bound} Let $d=1$ and $s>1/4$. Let $\mu$ be a minimizer of $E_{s,T}$ and assume that $\mu_T$ is upper $\alpha$-Ahlfors regular for some $\alpha>1-2s$. Then, for any $\eps\in(0,\T)$ and 
    \begin{equation}
        \label{e:betareg}
        \beta<\beta_{\textup{reg}}(s, \alpha)=  \frac{2s-1+\alpha}{2(1-s)+1-\alpha}
    \end{equation} there holds
     \[I(\mu, (T-\eps, T)) \lesssim_{\mu,\beta} \eps^{\beta}. \]
\end{thm} 
 To prove Theorem \ref{thm:lower_local_bound}, we first establish in Section \ref{subsec:interpolatedsobolev} some properties of upper $\alpha$-Ahlfors regular measures. In particular we prove that this class of measures embeds  into some negative Sobolev space and is stable with respect to McCann interpolation, see Lemma \ref{lem:reg_displacement}. This  result, which is a generalization of $L^\infty$ bounds on interpolated measures (see \cite[Theorem 8.7 \& Corollary 19.5]{Villani}), could be of independent interest.  Then, we build a competitor by taking a minimizer $\mu$ and shrinking by a factor its branches in the interval $[T-\eps,T]$. This keeps the perimeter part unchanged and reduces the kinetic energy. Since we assume $d=1$, the transport plan induced by $\mu$ is the optimal one and the measure that this new competitor irrigates is the McCann interpolant of $\mu_{T-\eps}$ and $\mu_T$. Thus, we can estimate the boundary energy by the Wasserstein distance (Lemma \ref{lem:H-s}) using interpolation inequalities and the additional regularity properties of upper $\alpha$-Ahlfors regular measures. \\
Let us now get to the proof of \eqref{e:alphamore}. In this case we have the following local scaling law.
\begin{thm} \label{thm:local_energy_lower}
    Let $d=1$ and $s>\frac{1}{4}$.  Let $\mu$ be a minimizer of $E_{s,T}$ and assume that  $\mu_T$ is lower $\alpha$-Ahlfors regular. Then, for any $\eps\in (0, T)$ and for  
    \begin{equation}
        \label{e:betacon}
        \beta_{\textup{con}}(\alpha) = \frac{2-\alpha}{2+\alpha}
    \end{equation} 
    there holds 
    \begin{equation}\label{e:upperboundth1.6}I(\mu, (T-\eps, T)) \lesssim_{M,\alpha} \eps^{\beta_{\textup{con}}}. \end{equation}
\end{thm}
The proof of Theorem \ref{thm:local_energy_lower} is based on the construction of a competitor  in $[T-\eps,T]$. To this extent we notice that lower $\alpha$-Ahlfors regular measure are supported on a $\alpha$ dimensional set. Moreover, if $d=1$, minimizers enjoy the so-called "cone property" (see Section \ref{subsec:1D}). This property guarantees a localization property of the branched transport plan. This allows us to make a local construction in the spirit of \cite[Theorem 1.3]{dephilippis2023energy} (see Lemma \ref{lem:single_branch}), from which Theorem \ref{thm:local_energy_lower} follows. \\

The paper is organized as follows: 
In Section \ref{sec:preliminaries} we recall some facts about Wasserstein distance and Sobolev spaces and we introduce the energy functional. Afterwards, we recall some already known properties of the functional and we highlight some properties which are specific to $d=1$. We then prove some global upper bounds. In Section \ref{sec:dim} we introduce the concepts of dimension of a measure and of upper and lower $\alpha$-Ahlfors regular measures. We then introduce the results of \cite{dephilippis2023energy} which we generalize to any $d\in \N$ and $s\in (0,1)$. Afterwards, assuming Theorems \ref{thm:lower_local_bound} and \ref{thm:local_energy_lower}, we prove Theorem \ref{thm:Malphareg}. Section \ref{sec:regular} is devoted to the study of upper $\alpha$-Ahlfors regular measures, followed by the proof of Theorem \ref{thm:lower_local_bound}. Finally in Section \ref{sec:concentrated} we adapt the construction of \cite{dephilippis2023energy} to prove Theorem \ref{thm:local_energy_lower}.

\section{Preliminaries}
\label{sec:preliminaries}
In this section, we introduce the mathematical model and the necessary notation. We start by recalling  basic facts about optimal transport and the Wasserstein distance.  Afterward, we introduce the internal energy and recall some results of \cite{CGOS}, highlighting results which hold only for $d=1$. Then, we recall some properties of fractional Sobolev spaces. Finally we introduce the full energy functional and derive global scaling laws.

\subsection{Notation} 
We use the symbols $\lesssim$ and $\gtrsim$ to indicate that there exists $C>0$ depending only on $d$ and $s$ such that  $a \leq C b$ or $a\geq C b$ respectively. Additionally, we denote by $ a \sim b$ a relation that satisfies $a\lesssim b$ and $a\gtrsim b$. We use $A\ll B$ as a hypothesis. It means that there exists an $\eps>0$ such that for $A\leq \eps B$ the conclusion holds.   For $x\in \R^d$,  we set $\vert x\vert_{per} = \min\{ \vert x - z\vert \colon z\in \Z^d  \}$ for the  distance on the  torus $\Tor^d= (\R/\Z)^d$.  We write $B_r(x)$ for a ball with respect to the Euclidean distance. We also  write $Q_r(x)$ for a square of side $r$ centered in $x$, parallel to the usual basis of $\R^d$. For $z\in \R^{d+1}$ we write $z=(x,t)\in \R^d \times \R$. The $\nabla$ operator  is always meant with respect to the variable $x$ and $\partial_t$ with respect to $t$. For $T>0$, we set $\Tor^d_T= \Tor^d \times (-T,T)$. We denote by $\mathcal{B}(\Tor^d)$ the Borel $\sigma$-algebra of $\Tor^d$. Let $A\subseteq \mathcal{\R}^d$ be a Borel set. We denote by $\mathcal{M}(A)$ the set of finite measures on $A$, by $\mathcal{M}^+(A)$ the set of positive measures, and by $\mathcal{P}(A)$ the set of probability measures on $A$. If $f$ is a function and $\mu$ a measure, we denote $f\#\mu$ the pushforward of $\mu$ by $f$.

Let $\rho_1\in C^{\infty}_c(B_1)$ denote a radially symmetric function such that $\rho_1\geq 0$ and  $\int_{\R^d}\rho_1\,dx=1$. Set $\rho_\eps(x)=\eps^{-d}\rho_1(x/\eps)$ so that for some $\eps_0>0$ small enough
 $\{ \rho_\eps \colon \eps \in (0,\eps_0)\}$ is a family of convolution kernels (extended periodically) that satisfy
\begin{equation}
    \label{prop:reg_conv_kernel}
    \rho_\eps \in C^{\infty}(\Tor^d), \quad \supp(\rho_\eps)\subseteq B_{\eps}, \quad \rho_\eps \geq 0, \quad \textup{ and }\int_{\Tor^d} \rho_\eps \,d x = 1 \quad \textup{ for all } \eps \in (0,\eps_0).
\end{equation}
For $\mu \in \mathcal{P}(\Tor^d)$, we set
\begin{equation}
    \label{def:conv_measure}
    (\rho_\eps\ast \mu)(x) = \int_{\Tor^d} \rho_\eps(x-y) d \mu(y) \quad \textup{ for } x\in \Tor^d.
\end{equation}

\subsection{Optimal transport}
In this section, we recall some basic facts and definitions from optimal transport. We refer to \cite{AGS,santambrogio2015optimal} for a broader overview of optimal transport. \\
If $\mu_1,\mu_2\in \mathcal{M}^+(\R^d)$ with $\mu_1(\R^d)=\mu_2(\R^d)$, we define 
\begin{equation}
\label{e:wasserstein}
    W^2(\mu_1,\mu_2)=\min\left\{\int_{\R^d\times \R^d}\abs{x-y}^2\,d\pi \colon \pi \in \mathcal{M}^+(\R^d\times \R^d), \ \pi_i=\mu_i \textup{ for } i=1,2\right\}.
\end{equation}
Analogously, if $\mu_1,\mu_2\in \mathcal{M}^+(\Tor^d)$ with $\mu_1(\Tor^d)=\mu_2(\Tor^d)$, we define 
\begin{equation}
\label{e:wassersteinperiodic}
    \Wper^2(\mu_1,\mu_2)=\min\left\{\int_{\Tor^d\times \Tor^d}\abs{x-y}_{\text{per}}^2\,d\pi \colon \pi \in \mathcal{M}^+(\Tor^d\times \Tor^d), \ \pi_i=\mu_i \textup{ for } i=1,2\right\}.
\end{equation}
 Let us recall the Benamou-Brenier formula: for $\mathbb{X}\in \{\Tor^d,\R^d\}$ and  $a,b\in\mathbb{R}$ with $a<b$ we have
\begin{equation}
\label{e:benamoubrenier}
    \frac{1}{b-a}W_*^2(\mu_1,\mu_2)=\min_{(\mu,m)}\left\{\int_{\X\times (a,b)}\abs{\frac{dm}{d\mu}}^2d\mu d t \colon  m\ll \mu, \ \partial_t\mu +\nabla\cdot m=0, \ \mu_a=\mu_1, \ \mu_b=\mu_2 \right\},
\end{equation}
where $W_*=W$ if $\X=\R^d$ and $W_*=W_{\textup{per}}$ if $\X=\Tor^d$.
Given $\mu_a,\mu_b\in \mathcal{M}^+(\R^d)$ with $\mu_a(\R^d)=\mu_b(\R^d)$ and $a\leq t\leq b$, we define the interpolated measure (McCann's interpolant) as 
\begin{equation}
\label{e:mccann}
\int_{\R^d}\psi d\mu^t=\int_{\R^d \times \R^d}\psi\left(\frac{b-t}{b-a}x+\frac{t-a}{b-a}y\right)d\pi
\end{equation}
 where $\pi$ is an optimal transport plan for $W(\mu_a,\mu_b)$.
This is the same as saying $\mu^t=F_t\#\pi$, where $F_t:\R^d\times\R^d \rightarrow\R$ is defined as \begin{equation*}
    F_t(x,y)=\frac{b-t}{b-a}x+\frac{t-a}{b-a}y.
 \end{equation*}
 By \cite[Theorem 5.27]{santambrogio2015optimal}, for $a\leq s\leq t\leq b$
 \begin{equation}
 \label{e:geodesic}
     W(\mu^t,\mu^s)=\left(\frac{t-s}{b-a}\right)W(\mu_a,\mu_b).
 \end{equation}
\subsection{The internal energy functional}
\label{subsec:energydef} 
Following the notation from \cite{dephilippis2023energy}, let us define the internal energy.
\begin{definition} 
    We denote by $\mathcal{A}(\Tor^d, (a,b))$ for $-T\leq a <b \leq T$ the set of measures $\mu\in \mathcal{M}^+(\Tor^d\times(a,b))$, $m\in\mathcal{M}(\Tor^d\times(a,b),\R^d)$, with $m\ll\mu$ satisfying in the sense of distributions the continuity equation 
    \begin{equation*}
        \label{e:continuityeq}
        \partial_t\mu+\nabla\cdot m=0 \quad \text{ in } \Tor^d \times (a,b),
    \end{equation*}
    and such that $\mu=\mu_t\otimes dt$ where, for almost every $t\in(a,b)$ there exist $N_t\in\N\cup \{+\infty\} $ such that it holds $\mu_t=\sum_{i=1}^{N_t}\varphi_i\delta_{X_i}$, for some $\varphi_i>0$ and $X_i \in \Tor^d$ with $i=1\dots N_t$. We denote by $\mathcal{A}^*(\Tor^d, (a,b))=\{\mu: \exists m ,\, (\mu, m)\in \mathcal{A}(\Tor^d, (a,b))\}$ the set of admissible $\mu$. For $\Tor^d_T = \Tor^d\times (-T, T)$, we use the notation $\mathcal{A}_T = \mathcal{A}(\Tor^d, (-T, T))$ and $\mathcal{A}_T^*= \mathcal{A}^*(\Tor^d, (-T, T)).$
    We define the functional $I:\mathcal{A}(\Tor^d, (a, b))\rightarrow[0,+\infty]$ by 
    \begin{equation*}
        I(\mu,m,(a,b))=\int_{a}^b\sum_{x\in \Tor^d}(\mu_t(\{x\}))^\frac{d-1}{d}d t + \int_{\Tor^d\times(a,b)}\abs{\frac{dm}{d\mu}}^2d \mu.
    \end{equation*}
    We also define (with abuse of notation) the functional $I:\mathcal{A}^*(\Tor^d, (a,b))\rightarrow[0,+\infty]$
    \begin{equation*}
        I(\mu,(a,b))=\min_m\{I(\mu,m,(a,b))\colon (\mu,m)\in \mathcal{A}_{\T} \}.
    \end{equation*}
    We use the notation $I(\mu)=I(\mu,(-T,T))$.
    We then set for $-\T\leq a< t<b\leq \T$
    \begin{equation*}
        \dot{P}(\mu,t)=\sum_{x\in\Tor^d}(\mu_t(\{x\}))^\frac{d-1}{d}, \quad \dot{E}_{\text{kin}}(\mu,t)=\int_{\Tor^d}\abs{\frac{dm_t}{d\mu_t}}^2d \mu_t
    \end{equation*}
    and
        \begin{equation*}
        P(\mu,(a,b))=\int_{a}^b\dot{P}(\mu,t)dt, \quad E_{\textup{kin}}(\mu,(a,b))=\int_{a}^b\dot{E}_{\text{kin}}(\mu,t)dt,
    \end{equation*}
    so that $I(\mu,(a,b))=P(\mu,(a,b))+E_{\textup{kin}}(\mu,(a,b))$.
\end{definition}

Notice that for $-\T\leq a <b\leq \T$ and $\mu\in \mathcal{A}^*_\T$, by \eqref{e:benamoubrenier} we have 
    \begin{equation}
        \label{e:ecinestimate}
        E_{\textup{kin}}(\mu,(a,b))\geq\frac{1}{b-a} \Wper^2(\mu_a,\mu_b).
    \end{equation}
Given two measures $\Bar{\mu}_\pm\in \mathcal{M}^+(\Tor^d)$ such that $\Bar{\mu}_+(\Tor^d)=\Bar{\mu}_-(\Tor^d)$
 consider the problem
\begin{equation}
\label{e:problembc}
    \inf \{I(\mu)\colon \mu \in \mathcal{A}^*_T, \ \mu_{\pm T}=\Bar{\mu}_\pm\}.
\end{equation}
By \eqref{e:ecinestimate}, if $I(\mu)<\infty$ then for $-\T< a<b< \T$
\begin{equation*}
    \Wper(\mu_a,\mu_b)\leq I(\mu)^\frac{1}{2} (b-a)^\frac{1}{2}.
\end{equation*}
This means that $t\mapsto\mu_t$ is a $\frac{1}{2}$-H\"older continuous curve in the Wasserstein space, hence $\mu_T$ and $\mu_{-T}$ are well-defined.\\
Let us now recall some results of \cite{CGOS,dephilippis2023energy} for this functional. First, arguing as in \cite[Proposition 5.5]{CGOS} we have existence of  minimizers (in \cite{CGOS} the discussion is carried for $d=2$, but the argument extends to any dimension). Moreover, when $\Bar{\mu}_+=\Bar{\mu}_-$, by a symmetrization argument we can concentrate on minimizers which are symmetric with respect to $t=0$. Let us recall the notion of subsystems, see \cite[Proposition 2.6]{dephilippis2023energy} and \cite[Proposition 5.7]{CGOS}.
\begin{proposition}[Subsystems]\label{prop:subsys}
Let $t\in[-T,T]$ and $(\mu,m)\in \mathcal{A}_T$ with $I(\mu,m)<\infty$. Set $v=dm/d\mu$. Then, for every measure $\sigma\leq \mu_t $ there exists a measure $\Tilde{\mu}\in\mathcal{A}^*_T$ such that
\begin{enumerate}
    \item $\mu'=\mu-\Tilde{\mu}\in\mathcal{M}^+(\Tor^d_T)$;
    \item $\Tilde{\mu}_t=\sigma$;
    \item $\Tilde{\mu}$ satisfies the continuity equation 
    $        \partial_t\Tilde{\mu}+\nabla\cdot (v\Tilde{\mu})=0 $.
\end{enumerate}

In particular $(\Tilde{\mu},v\Tilde{\mu})\in \A_T$ and $I(\Tilde{\mu},v\Tilde{\mu})\leq I(\mu,m)$. Similarly $(\mu',v\mu')\in \A_T$ and $I(\mu',v\mu')\leq I(\mu,m)$.
\end{proposition}

If $\mu_t=\sum_i \varphi_i\delta_{X_i}$ and $\sigma=\varphi_i\delta_{X_i}$ for some $i$,
we call $\mu_+^{t,i}=\Tilde{\mu}\restr(t,T)\times \Tor^d$ the forward subsystem emanating from $X_i$ and $\mu_-^{t,i}=\Tilde{\mu}\restr(-T,t)\times \Tor^d$ the backward subsystem emanating from $X_i$.

Let us now recall the no-loop property of minimizers, see \cite[Lemma 5.8]{CGOS}.
\begin{lemma}[No-loop condition]
\label{prop:noloop}
Let $\mu$ be a minimizer of \eqref{e:problembc} and $t\in(-T,T)$. Let $X_1,X_2\in \Tor^d$ with $\mu_t(X_i)\neq 0$ for $i=1,2$ and consider the forward and backward subsystems $\mu_\pm^{t,i}$ emanating from $X_i$. If there exist $-T\leq t_-<t<t_+\leq T$ such that for $*\in\{\pm\}$, $\mu_*^{t_*,1}$ and $\mu_*^{t_*,2}$ are not mutually singular, then  $X_1=X_2$.

\end{lemma}
 Let us also recall the following representation formula, see \cite[Proposition 5.10 and Proposition 5.11]{CGOS}.
\begin{proposition}
    Let $\Bar{\mu}\in \mathcal{M}^+(\Tor^d)$ and $\mu$ be a symmetric minimizer of \eqref{e:problembc} with $\mu_\pm=\Bar{\mu}$. Then, for every $T'<T$, $\mu$ is a finite graph in $\Tor^d_{T'}$ without loops. In particular we have

        \begin{equation}
        \label{e:energyrepresentation}
I(\mu) = \int_{-T}^T\sum_{i}\left[\varphi_i^{\frac{d-1}{d}}(t)+ \varphi_i(t) \vert \dot{X}_i(t)\vert^2 \right]\, dt
\end{equation}
where the sum is locally finite. Moreover, between two branching points $t\mapsto\varphi_i(t)$ is constant while $t\mapsto X_i(t)$ is affine.
\end{proposition}
In the rest of the paper we will always use the representation formula \eqref{e:energyrepresentation} for $I$. 
Let us now introduce a Lagrangian reformulation of the problem. The proof of the following proposition is essentially the same as \cite[Theorem 4.3]{dephilippis2023energy} so we omit it. We just remark that it is essential that the no-loop property holds so that we can parametrize the curves induced by a minimizer $\mu\in\mathcal{A}^*_T$ by their initial point. This avoids the introduction of traffic plans like in \cite{pegonfractal2018,landscape}. Let
\begin{equation*}
    C_T=\{X:\Tor^d\times[-T,T]\mapsto \Tor^d\,: X(x,-T)=X(x,T)=x \textup{ and for a.e.}\ x \ t\mapsto X(x,t) \textup{ is AC}\},
\end{equation*}
where AC denotes absolutely continuous curves. Let $X\in C_T$, $\overline{\mu}\in \mathcal{P}(\Tor^d)$, $y\in \Tor^d$ and $t\in (-T,T)$. We define the multiplicity function as
\begin{equation*}
    \varphi_{X,\overline{\mu}}(y,t)=\overline{\mu}(\{x\in \supp \overline{\mu}:\ X(x,t)=y\}).
\end{equation*}
We then set
\begin{equation*}
    \mathcal{L}(X,\overline{\mu})  =\int_{\Tor^d_T}\lt( \varphi_{X,\overline{\mu}}(X(x,t),t)^{-\frac{1}{d}}+\abs{\partial_t X(x,t)}^2\rt)\,d\overline{\mu}(x)dt.
\end{equation*}
\begin{proposition}
\label{prop:lagrangian}
    For every $T>0$ and $\overline{\mu}\in \mathcal{P}(\Tor^d)$
    \begin{equation*}
        \min_{\mu\in \A^*_T}\{I(\mu,(-T,T)),\mu_{\pm T}=\overline{\mu}\}=\min_{X\in C_T}\mathcal{L}(X,\overline{\mu}).
    \end{equation*}
    Moreover, for every $X$ minimizing $\mathcal{L}$, $\mu_t=X(\cdot,t)\#\overline{\mu}$ is a minimizer for $I$. Viceversa, we can associate to any minimizer $\mu\in  \A^*_T$ a collection of curves $X\in C_T$ (defined for $\overline{\mu}$-a.e. $x\in\Tor^d$) such that $X$ is a minimizer for $\mathcal{L}$.
\end{proposition}
Let us finally recall the equipartition of energy \cite[Proposition 2.11]{dephilippis2023energy}.
\begin{proposition}[Equipartition of Energy]

    \label{prop:equipartition}
    Let $\Bar{\mu}\in \mathcal{M}^+(\Tor^d)$ and consider a symmetric minimizer of \eqref{e:problembc} with $\mu_\pm=\Bar{\mu}$. 
    There exists $\Lambda\in \mathbb{R}$ with
    \begin{equation}
        \abs{\Lambda}\lesssim \frac{I(\mu)}{T}
    \end{equation}
    such that for every $t\in(-T,T)$
    \begin{equation}
       \dot{P}(\mu,t)=\dot{E}_{\text{kin}}(\mu,t)+\Lambda.
    \end{equation}
\end{proposition}

\subsection{The case $d=1$}
\label{subsec:1D}
When $d=1$, minimizers of \eqref{e:problembc} satisfy additional properties. A first simple but useful implication of the no-loop condition for $d=1$ is that the coupling between $\mu_T$ and $\mu_t$ induced by $\mu$ is monotone, see \cite[Lemma 4.1]{dephilippis2023energy}. 
\begin{lemma}
    \label{lem:monotonecoupling}Let $\Bar{\mu}\in \mathcal{M}^+(\Tor)$ and $\mu$ be a symmetric minimizer of \eqref{e:problembc} with $\mu_\pm=\Bar{\mu}$. For every $t\in[0,T)$, if  $\mu_t=\sum_i \varphi_i\delta_{X_i}$ and if we denote by $\mu^{t,i}=\mu^{t,i}_+$ the forward subsystems emanating from $X_i$, then considering $\mu_t$ as a periodic measure on $\R$ we have
    \begin{equation}
        X_i<X_j \quad \implies \quad \sup\supp \mu_T^{t,i} \leq \inf\supp \mu_T^{t,j}.
    \end{equation}
    Moreover for $i\neq j$, \begin{equation*}
        \mu_T^{t,i}(\{x\})\mu_T^{t,j}(\{x\})=0 \qquad \forall x \in \Tor .
    \end{equation*}
    As a consequence, the map $X(\cdot,t)$ given by Proposition \ref{prop:lagrangian} (considered as a map on $\R$) is monotone.
\end{lemma}
 We now recall from \cite{dephilippis2023energy} the cone property. Let us consider the parametric representation $X$ of a minimizer $\mu$ of \eqref{def:total_energy} given by Proposition \ref{prop:lagrangian}, such that for any $x\in\Tor$, $X(x,T)=X(x,-T)=x$ and for any $t\in (-\T,\T)$ it holds $X(\cdot,t)\sharp\mu_T=\mu_t$. Let us define the irrigated set of a given point $\bar{x}$ at time $t\ge 0$ as
\begin{equation*}
    J(\Bar{x},t)=\{x\in\Tor\, |\,  X(x,t)=\Bar{x}\}.
\end{equation*}
 Let   $C(\Bar{x},t)$ be the cone defined by the convex hull in space and time of the point $(\Bar{x},t)$ and the set $J(\Bar{x},t)\times\{T\}$. The following was proved in \cite[Proposition 4.7]{dephilippis2023energy}.
\begin{proposition}[Cone property]\label{prop:cone}
    Let $x\in J(\Bar{x},t)$. Then $X(x,s)\in C(\Bar{x},t)$ for any $0<s<t$.
\end{proposition}
We now observe, as in \cite{minimizers2d} and \cite[Proposition 4.10 \& Remark 4.11]{dephilippis2023energy}, that for $t=0$, if $\mu_0=\sum_i \varphi_i \delta_{X_i}$ then $X_i$ must coincide with the barycenter of the corresponding irrigated measure $\mu_T^i= (\mu^{0,i}_+)_T$.
\begin{lemma}\label{lem:barycenter}
Let $\bar \mu\in \M^+(\Tor)$ and  $\mu$ be a symmetric minimizer of \eqref{e:problembc}.  Extending $\mu$ by periodicity on $\R$, letting $\mu_0=\sum_i \varphi_i \delta_{X_i}$ and  $\mu_T^i= (\mu^{0,i}_+)_T$, we have for every $i$,
\begin{equation}\label{bary}
 \varphi_i X_i= \int_{\R} y d\mu_T^i.
\end{equation}
\end{lemma}
\begin{remark}
\label{rem:periodicnonperiodic}
    As a consequence of Lemma \ref{lem:barycenter} and the cone property, up to choosing appropriately the origin, a periodic minimizer $\mu$ of $I$ can be identified with a minimizer $\Tilde{\mu}$ of $I$ on $[0,1]\times[-T,T]$ without boundary conditions. This in particular implies that \eqref{e:ecinestimate} holds also for the Euclidean Wasserstein distance \eqref{e:wasserstein}. Moreover, if we now consider a parametric representation $\tilde{X}$ of $\tilde{\mu}$ as given by Proposition \ref{prop:lagrangian}, by Lemma \ref{lem:monotonecoupling} we have that $\Tilde{X}(\cdot,t)$ is a monotone map on $\supp\Tilde{\mu}_T$. Since $d=1$, this implies that for all $t\in(0,T) $ $ \Tilde{X}(\cdot,t)$ is an optimal transport map for $W(\tilde{\mu}_T,\tilde{\mu}_t)$. In the rest of the paper we will omit distinguishing between $\mu$ and $\Tilde{\mu}$.
\end{remark}

\begin{figure}[h]
    \centering\includegraphics[scale = 1]{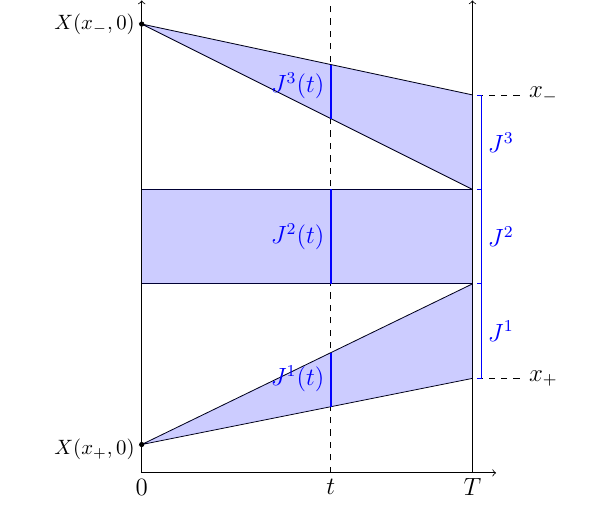}
    \caption{ Representation of the typical case for the $3$ intervals $J^1(t)$, $J^2(t)$ and $J^3(t)$ for $t\in [0,T]$. }
    \label{fig:3intervals}
\end{figure}

\begin{proof}
Let $X$ be a parametric representation of $\mu$ as given by Proposition \ref{prop:lagrangian}. As in \cite[Proposition 4.10]{dephilippis2023energy} we define 
 \[
  \widehat{X}(x,t)=X(x,t)- \frac{T-t}{T} \left(X_i-\frac{1}{\varphi_i}\int_{\R}y\,d\mu_T^i\right) \qquad \textrm{for } t\in[0,T] \textrm{ and } x\in \supp \mu_T^i.
 \]
 We then extend $\widehat{X}$ to $[-T,T]$ by symmetry. Notice that by the cone property this construction leaves the perimeter unchanged. Hence by a direct computation
 
\[
 \mathcal{L}(\widehat{X},\bar \mu)=\mathcal{L}(X,\bar \mu)- \frac{\varphi_i}{T}\lt| X_i -\frac{1}{\varphi_i}\int_{\R} y d\mu_T^i\rt|^2.
\]
This concludes the proof.
\end{proof}
We now want to show that given an interval $J=(x^+,x^-)\subset \Tor$, it can be irrigated by (or rather it "splits in") at most 3 disjoint intervals. 
We let  
\begin{align*}
    &J^1=\{x\in J  \, | \, \ \exists t\in (0,\T)\, \textrm{ s.t. }  X(x,t)=X(x^+,t)\},\\
     &J^2=\{x\in J  \,| \,  \  X(x,t)\neq X(x^\pm,t)\,  \, \forall t\in (0,T) \},\\
    &J^3=\{x\in J  \, | \, \ \exists t\in (0,\T)\, \textrm{ s.t. }  X(x,t)=X(x^-,t) \}.
\end{align*}
Clearly $J^2$ does not intersect with $J^1$ and $J^3$ , and may be empty. By the no-loop condition, if $J^1$ and $J^3$ intersect, then they must coincide and in that case $J^2$ must be empty. In the case where the three sets are disjoined and non-empty, again by the no loop condition they must satisfy $\sup J^1\leq \inf J^2$ and $\sup J^2\leq \inf J^3$, see Figure \ref{fig:3intervals}. Notice that if $\Bar{x}$ is such that $J(\Bar{x},0)\subset J^2$, by Lemma \ref{lem:barycenter}, we must have $\Bar{x}\in J^2$ and thus by Proposition \ref{prop:cone} for $t=0$, also $ C(\Bar{x},0)\subset J^2\times[0,T]$. This implies in particular:
\begin{proposition}
    For every $x\in J^2$ and every $t\in [-T,T]$, $X(x,t)\in J^2$.
\end{proposition}
%
%
Moreover, by the cone property, we have that for all $x\in J^1$, $X(x,t)\in C(X(x^+,T),T)$ and for all $x\in J^3$, $X(x,t)\in C(X(x^-,T),T)$. This gives us the following corollary, see Figure \ref{fig:3intervals}.
\begin{corollary}
\label{cor:3intervals}
    Let $J\subset  \Tor$ be an interval. Then, for any $t\in [0,T]$ there exists at most three disjoint intervals $J^i(t)$, $i=1,2,3$, with $|J^i(t)|\leq |J|$  and such that for $\mu_T$-a.e.  $ x \in J$ we have $X(x,t)\in J^i(t)$ for some  $i=1,2,3$.
\end{corollary}

\subsection{Negative Sobolev spaces}
Let us recall some facts regarding homogeneous negative Sobolev spaces, see \cite{bahourichemindanchin,triebel}. For $s\in \R$, we define $\Hdots(\Tor^d)$ as the set of all  measures $\sigma\in \mathcal{M}(\Tor^d)$ such that
\begin{equation}
    \label{def:hom_measures_negative_sobolev}
    \Vert \sigma\Vert_{\Hdots}^2 = \sum_{k\in \Z^d} \vert k\vert^{-2s} \vert \widehat{\sigma}_k\vert^2 <\infty,
\end{equation}
where $\widehat{\sigma}_k$ is the Fourier coefficient for $k\in \Z^d$ given by 
\begin{equation}
    \label{def:fourier}
    \widehat{\sigma}_k = \int_{\Tor^d} \exp(
    2\pi i(k\cdot x))d \sigma(x). 
\end{equation}
Notice that with this notation, $\dot{H}^{-0}(\Tor^d)=L^2(\Tor^d)$ and that for $s<0$ \eqref{def:hom_measures_negative_sobolev} implies $\sigma(\Tor^d)=0$. We recall that $\dot{H}^{-s}(\Tor^d)=(\dot{H}^s(\Tor^d))^*$ and that for $s\in (0,1)$, $\dot{H}^s(\Tor^d)$  coincides  with the completion of the set of  $C^\infty(\Tor^d)$ functions with zero mean with respect to the norm $\norm{\psi}_{\dot{H}^s(\Tor^d)}^2$. Moreover, the exist a constant $C(d,s)$ such that
\begin{equation}
\label{e:sobolevseminorm}
    \norm{\psi}_{\dot{H}^s(\Tor^d)}^2=C(d,s)\int_{\Tor^d}\int_{\Tor^d}\frac{\abs{\psi(x)-\psi(y)}^2}{\abs{x-y}_{\textup{per}}^{d+2s}}\,dxdy.
\end{equation}
We will also need the space $\dot{H}^s(\Omega)$ for a bounded Lipschitz domain $\Omega\subset \R^d$, defined as the completion of $C^\infty(\Omega)$ with respect to the norm
\begin{equation}
\label{e:sobolevseminorm2}
    \norm{\psi}_{\dot{H}^s(\Omega)}^2=C(d,s)\int_\Omega\int_\Omega\frac{\abs{\psi(x)-\psi(y)}^2}{\abs{x-y}^{d+2s}}\,dxdy.
\end{equation}
Let us present some characterizations of the negative Sobolev norm.
\begin{proposition}Let $s\in (0,1)$ and $(\rho_\eps)_{\eps>0}$ be a family of regular symmetric convolution kernels  with $\eps\in(0,\eps_0)$ as in \eqref{prop:reg_conv_kernel}. For all $\sigma\in \mathcal{M}(\Tor^d)$ with $\sigma(\Tor^d)=0$, and $\gamma>0$ there holds
\begin{equation}
\label{e:H-scharacterization1}
    \Vert \sigma\Vert_{\Hdots}^2 \lesssim_\gamma \int_0^1 \eta^{2s+\gamma} \left( \sum_{\vert k\vert \leq \eta^{-1}} \vert k \vert^\gamma \vert \widehat{\sigma}_k\vert^2 \right) \frac{d \eta}{\eta}
\end{equation}
and 
\begin{equation}
\label{e:H-scharacterization2}
    \Vert \sigma\Vert^2_{\dot{H}^{-s}} \lesssim \int_0^{\eps_0}\eps^{2s} \int_{\Tor^d} \vert \rho_\eps \ast \sigma\vert^2 d x \frac{d\eps}{\eps}.
\end{equation}
\end{proposition}
\begin{proof}
    For a proof of the first estimate with $\gamma=1$, we refer to \cite[Proposition 2.13]{dephilippis2023energy}. The proof is essentially the same for a generic $\gamma>0$. Let us prove the second estimate, see also \cite[Theorem 2.7]{GolHues} for a variant. First, notice that by Plancherel's identity and properties of Fourier coefficients we have
    \begin{equation*}
         \int_{\Tor^d} \vert \rho_\eps \ast \sigma\vert^2 d x=\sum_{k\in \mathbb{Z}^d\backslash \{0\}}\abs{\widehat{\sigma}_k}^2|\widehat{(\rho_\eps)}_k|^2.
    \end{equation*}
    Thus, to conclude we just need to show that 
    \begin{equation*}
        \int_0^{\eps_0}\eps^{2s}|\widehat{(\rho_\eps)}_k|^2
\,\frac{d\eps}{\eps}\gtrsim \abs{k}^{-2s}.
    \end{equation*}
    With a change of variable we have (recall \eqref{prop:reg_conv_kernel})
    \begin{multline*}
         \widehat{(\rho_\eps)}_k=\int_{\Tor^d}\rho_\eps(x) \exp(
    2\pi i(k\cdot x))d x=\int_{B_{\eps}}\frac{1}{\eps^d}\rho_1(\frac{x}{\eps}) \exp(
    2\pi i(k\cdot x))d x\\=\int_{ \R^d}\rho_1(x) \exp(
    2\pi i(\eps k\cdot x))d x=F(\eps k).
    \end{multline*}
     Let us notice that, since $\rho_1$ is radially symmetric,  $F(v)=\tilde{F}(\abs{v})$ for  $v\in\R^d$.
    Finally we have
    \[
        \int_0^{\eps_0}\eps^{2s}\abs{F(\eps k)}^2\frac{d\eps}{\eps}=\abs{k}^{-2s}\int_{0}^{\abs{k}\eps_0}t^{2s}\abs{F(tk/\abs{k})}^2\frac{dt}{t} \geq \abs{k}^{-2s}\int_{0}^{\eps_0}t^{2s}\abs{\tilde{F}(t)}^2\frac{dt}{t}.
    \]
    Since $\tilde{F}(0)=1$ and $\tilde{F}$ is continuous on the integration domain, the last integral is bounded from below.

\end{proof}
We will use the following rather standard result.
\begin{proposition}
\label{prop:kernelestimate}
    Let $s\in (0,1)$, $\gamma \in (0,s]$, and $(\rho_\eps)_{\eps>0}$ be a family of regular convolution kernels with $\eps\in(0,\eps_0)$ as in \eqref{prop:reg_conv_kernel}. For every $\sigma \in \mathcal{M}(\Tor^d)$ with $\sigma(\Tor^d)=0$
    \[\Vert \sigma -\sigma \ast \rho_\eps \Vert_{\Hdots} \lesssim \eps^{s-\gamma} \Vert \sigma \Vert_{\dot{H}^{-\gamma}}.\]
\end{proposition}
\begin{proof}

we have for $k\neq 0$ and using $(\widehat{\sigma\ast \rho_\eps})_k=\widehat{\sigma}_k(\widehat{\rho}_\eps)_k$
\begin{equation}
 \label{est:interpolation_fractal}
        \Vert \sigma - \sigma \ast \rho_\eps \Vert_{\Hdots}^2 = \sum_{\vert k\vert \leq 1/\eps} \vert k\vert^{-2s} \vert \widehat{\sigma}_k\vert^2 \vert 1 - (\widehat{ \rho}_\eps)_k\vert^2 + \sum_{\vert k \vert > 1/\eps} \vert k\vert^{-2s}\vert \widehat{\sigma}_k\vert^2 \vert 1 - (\widehat{ \rho}_\eps)_k\vert^2 = T_1 + T_2.
\end{equation}
We estimate the two terms on the right-hand side of \eqref{est:interpolation_fractal} separately. We start with $T_2$. If  $\vert k\vert \geq 1/\eps$ then  $\vert k\vert^{-2(s-\gamma)} \leq \eps^{2(s-\gamma)}$. Moreover, for every $k\in \N$,
\begin{equation*}
    |(\widehat{ \rho}_\eps)_k|=\abs{\int_{\Tor^d}\exp(2\pi i k\cdot x)\rho_\eps(x)\,dx}\leq \int_{\Tor^d}\rho_\eps(x)\,dx=1.
\end{equation*}
Hence $T_2$ is bounded by
\begin{equation}
    \label{bound_T2}
    T_2 = \sum_{\vert k \vert > 1/\eps} \vert k\vert^{-2s} \vert \widehat{\sigma}_k\vert^2 \vert 1 - (\widehat{ \rho}_\eps)_k\vert^2 \lesssim  \eps^{2(s-\gamma)} \sum_{\vert k\vert > 1/\eps} \vert k\vert^{-2\gamma} \vert \widehat{\sigma}_k\vert^2 \lesssim \eps^{2(s-\gamma)} \Vert \sigma \Vert_{\dot{H}^{-\gamma}}^2.
\end{equation}
For $T_1$, we  recall that for $\abs{z}\lesssim 1$ we have $\abs{1-\exp(z)}\lesssim \abs{z}$. Hence for $\abs{k}<1/\eps$ we estimate 
\begin{multline*}
    \vert 1 - (\widehat{ \rho}_\eps)_k\vert=\abs{\int_{\Tor^d}(1-\exp(2\pi i k\cdot x))\rho_\eps(x)\,dx}\leq \int_{\Tor^d}\vert 1-\exp(2\pi i \eps k\cdot x)\vert\rho_1(x)\,dx \\ \lesssim
   \int_{\Tor^d}\vert \eps k\cdot x\vert\rho_1(x)\,dx \lesssim \eps \abs{k}.
\end{multline*}
This leads to the bound 
\begin{equation}
    T_1 = \sum_{\vert k\vert \leq 1/\eps} \vert k\vert^{-2s}  \vert \widehat{\sigma}_k\vert^2 \vert 1 - (\widehat{ \rho_\eps})_k\vert^2\leq \eps^2 \sum_{\vert k\vert \leq 1/\eps} \vert k\vert^{2 -2s + 2\gamma} \vert k\vert^{-2\gamma} \vert \widehat{\sigma}_k\vert^2 \leq  \eps^{2(s-\gamma)} \Vert \sigma\Vert_{\dot{H}^{-\gamma}}^2. \label{bound_T1}
\end{equation}
Combining \eqref{bound_T1} and \eqref{bound_T2} implies the desired inequality.
\end{proof}
Finally, for $s=1$, the $\dot{H}^{-1}$ norm morally behaves like the Wasserstein distance for measures. We recall the following inequality, which is a special case of \cite[Proposition 2.8]{LOEPER200668}, see also \cite{peyre}.
\begin{proposition}
    \label{prop:H-sW2}
    Let $\sigma_1,\sigma_2\in L^\infty(\Tor^d)\cap \mathcal{M}^+(\Tor^d)$. Then,
    \begin{equation}
        \label{e:H-sW2}
        \Vert\sigma_1-\sigma_2\Vert_{\dot{H}^{-1}}\leq \max\{\Vert\sigma_1\Vert_{L^\infty},\Vert\sigma_2\Vert_{L^\infty}\}^\frac{1}{2}W_{\textup{per}}(\sigma_1,\sigma_2).
    \end{equation}
\end{proposition}

\subsection{The main functional and global scaling laws}
In this section, we introduce the full energy functional and recall some  results on the global scaling law. For $s\in (0,1)$ and $\T>0$ we define
 \begin{equation}\label{def:total_energy}
    E_{s, \T}(\mu) = I(\mu) + \Vert \mu_{\pm\T} - 1\Vert_{\dot{H}^{-s}}^2,
\end{equation} 
where 
\begin{equation*}
    \Vert \mu_{\pm\T} - 1\Vert_{\dot{H}^{-s}}^2=\Vert \mu_{\T} - 1\Vert_{\dot{H}^{-s}}^2+\Vert \mu_{-\T} - 1\Vert_{\dot{H}^{-s}}^2 .
\end{equation*}
By a  symmetrization argument, we may restrict the problem \eqref{def:total_energy} to minimizers that are symmetric with respect to $t=0$.

Let us now turn to global upper bounds for the energy.
\begin{proposition}[Global upper bounds]
\label{pro:global_scaling}
Let $T\geq 1$. We have
\begin{equation}
\label{e:thicksample}
    \min\{ E_{s,T}(\mu)\colon \mu \in \mathcal{A}^*_\T\} \lesssim  T.
\end{equation}  
Let $0 < T \leq 1 $ and  $d>2s$. For $s>\frac{1}{4}$ we have
    \begin{equation}
    \label{e:scalinggloballocal}
        \min\{ E_{s,T}(\mu)\colon \mu \in \mathcal{A}^*_\T\} \lesssim T^{\frac{d+2s}{3d+2(1-s)}}    
        \end{equation} 
    and for $s\leq\frac{1}{4}$ we have 
    \begin{equation}
    \label{e:scalinglebesgue}
    \min\{ E_{s,T}(\mu)\colon \mu \in \mathcal{A}^*_\T\} \lesssim T^{\frac{1}{3}}.
      \end{equation}
Finally if  $0 < T \leq 1 $, $d<2s$ we have
\begin{equation}
\label{e:scalingdelta}
     \min\{ E_{s,T}(\mu)\colon \mu \in \mathcal{A}^*_\T\} \lesssim  T^{\frac{2s}{2s+1}}.
\end{equation}
\begin{remark}
 While we do not use these bounds directly, they give indication on the behavior of the functional in the different parameter regimes and generalize \cite[Proposition 3.1]{dephilippis2023energy}. In particular, it is interesting to compare these global bounds with the local bounds which we find in Section 3. The most interesting regime is \eqref{e:scalinggloballocal}. Notice that \eqref{e:scalinggloballocal} coincides with the conjectured optimal local scaling \eqref{e:critical}. Notice also that if $d=2$ and $s=1/2$, it coincides with the results of \cite{Choksi2004} and \cite[Theorem 4.1]{CoOtSe2016}. The regime \eqref{e:thicksample} corresponds to samples which are too thick so that the details of the branching are not relevant. The regime \eqref{e:scalinglebesgue} instead corresponds to strong boundary penalizations, which force the boundary measure to have full dimension, see Remark \ref{rem:dimension}. Indeed the scaling \eqref{e:scalinglebesgue} matches the local scaling given by Proposition \ref{prop:scaling13}, which is the scaling of a uniformly distributed measure. Finally, in the regime \eqref{e:scalingdelta} Dirac deltas have finite $H^{-s}$ norm. In this case we do not expect that global and local scaling laws coincide. We also remark that, like in \cite{CGOS,CoOtSe2016,cintiotto}, it should be possible to find matching lower bounds, but this goes beyond the scope of the paper.
\end{remark}

\end{proposition}
\begin{proof}
Consider $N$ equidistributed squares of side length $r>0$, which we denote by $Q_r(X_i)$ where for $i=1,\dots, N$ $X_i$ denotes the centers. Let
 \begin{equation*}
    \mu_T=\sum_i \frac{1}{Nr^d}\chi_{Q_r(X_i)}.
\end{equation*}
 Following the strategy in \cite{CoOtSe2016} and \cite[Proposition 5.2]{CGOS}, it is possible to construct $\mu\in \mathcal{A}^*_T$ such that 
\begin{equation}
\label{e:scalingconstruction}
    E_{s,T}(\mu)\lesssim TN^\frac{1}{d}+\frac{r^2}{T}+ \norm{\mu_T-1}_{\dot{H}^{-s}}^2.
\end{equation}
If $T\ges 1$ we can just chose $N=r=1$ which gives $\mu_T=1$ and
\begin{equation*}
     E_{s,T}(\mu)\lesssim T+\frac{1}{T}\lesssim T.
\end{equation*}
 From now on we assume that  $T\ll 1$.
Using Proposition \ref{prop:boundaryest} below we get for $d>2s$
\begin{equation*}
     E_{s,T}(\mu)\lesssim TN^\frac{1}{d}+\frac{r^2}{T}+ \frac{1}{Nr^{d-2s}}.
\end{equation*}
Optimizing in $r$ we find the optimal value 
\begin{equation}
\label{e:optimalr}
    r=T^{\frac{1}{d+2(1-s)}}N^{-\frac{1}{d+2(1-s)}}
\end{equation}
that gives 
\begin{equation*}
    E_{s,T}(\mu)\lesssim TN^\frac{1}{d}+T^{-\frac{d-2s}{d+2(1-s)}}N^{-\frac{2}{d+2(1-s)}}.
\end{equation*}
Optimizing now in $N$ we find 
\begin{equation}
\label{e:optimaln}
    N=T^{-\frac{2d(1+d-2s)}{3d+2(1-s)}}
\end{equation}
which gives the bound
\begin{equation*}
    E_{s,T}(\mu)\lesssim T^\frac{d+2s}{3d+2(1-s)}.
\end{equation*}
Notice that since $T\ll1$, the choice \eqref{e:optimaln} of $N$ yields $N\gg1$ which is compatible with the requirement $N\in \N$. Plugging this choice of $N$ in \eqref{e:optimalr} we get
\begin{equation}
\label{e:optimalr2}
    r=T^\frac{1+2d}{3d+2(1-s)},
\end{equation}
and comparing \eqref{e:optimaln} and \eqref{e:optimalr2} the condition $r\ll N^{-1/d}$ becomes 
\begin{equation*}
    s>\frac{1}{4}.
\end{equation*}
If instead $s\leq \frac{1}{4}$, we choose $r^d=\frac{1}{N}$ in \eqref{e:scalingconstruction}, so that the boundary term is zero. This implies
\begin{equation*}
    E_{s,T}(\mu)\lesssim TN^\frac{1}{d}+\frac{1}{TN^\frac{2}{d}}.
\end{equation*}
Optimizing in $N$ we find the desired bound. \\
Finally suppose that $d<2s$. This can happen only for $d=1$ since we assume that $s\in(0,1)$. In this case, Dirac deltas are in $H^{-s}(\Tor)$. Consider $\mu_T$ to be a sum of $N$ equispaced Dirac deltas with weight $1/N$ and let $\mu_t=\mu_T$ for all $t\in [-T,T]$. we find
\begin{equation*}
    E_{s,T}(\mu)=2TN+\norm{\mu_T-1}_{\dot{H}^{-s}}^2.
\end{equation*}
By Proposition \ref{prop:deltasinH-s} below, we can estimate the boundary term to find 
\begin{equation*}
    E_{s,T}(\mu)\lesssim TN+\frac{1}{N^{2s}}.
\end{equation*}
Optimizing in $N$ leads to the thesis. 
\end{proof}
Let us estimate the boundary term in the following propositions. The first one is a slight variation of \cite[Lemma 4.4]{CoOtSe2016}.
\begin{proposition}
\label{prop:boundaryest}
Let $d>2s$, $N\gg 1$ and $0<r\leq N^{-1/d}$. Consider $N$ equidistributed squares of side length $r$, which we denote by $Q_r(X_i)$ where $X_i$, $i=1,\dots, N$ are the centers and let
 \begin{equation*}
    \mu=\sum_i \frac{1}{Nr^d}\chi_{Q_r(X_i)}.
\end{equation*} Then
\begin{equation}
    \norm{\mu-1}_{\dot{H}^{-s}}^2\lesssim \frac{1}{Nr^{d-2s}}.
\end{equation}
\end{proposition}
\begin{proof}
Set $\sigma=\mu-1$ and let $\psi\in L^1(Q_{N^{-1/d}}(X_i))$. Denote by $\psi_0$ the mean of $\psi$ on $Q_{N^{-1/d}}(X_i)$.
 Since  $\sigma(Q_{N^{-1/d}}(X_i))=0$, we have
\begin{equation*}
\int_{Q_{N^{-1/d}}(X_i)}\psi\,d\sigma=\int_{Q_{N^{-1/d}}(X_i)}(\psi-\psi_0)\,d\sigma=\int_{Q_r(X_i)}\frac{1}{Nr^d}(\psi-\psi_0)\,dx\leq \frac{1}{Nr^d}\norm{\psi-\psi_0}_{L^1(Q_r(X_i))}.
\end{equation*}
By Sobolev embedding, for every $\ell>0$, if $Q_\ell$ is a cube of side length $\ell$, $\dot{H}^s(Q_\ell)\hookrightarrow L^{p^*}(Q_\ell)$ where $p^*=2d/(d-2s)$, see for instance \cite{sobolevguide}. Moreover, the norm of the embedding does not depend on $\ell$.
Hence, by H\"older's inequality,
\begin{multline*}
    \norm{\psi-\psi_0}_{L^1(Q_r(X_i))}\leq \norm{\psi-\psi_0}_{L^{p^*}(Q_r(X_i))}\abs{Q_r(X_i)}^{1-\frac{1}{p^*}}\\ =r^{\frac{d+2s}{2}}\norm{\psi-\psi_0}_{L^{p^*}(Q_r(X_i))}\le r^{\frac{d+2s}{2}}\norm{\psi-\psi_0}_{L^{p^*}(Q_{N^{-1/d}}(X_i))}\lesssim r^{\frac{d+2s}{2}}\norm{\psi}_{\dot{H}^s(Q_{N^{-1/d}}(X_i))}.
\end{multline*}
Thus,
\begin{equation*}
    \int_{Q_{N^{-1/d}}(X_i)}\psi\,d\sigma\lesssim \frac{1}{Nr^{\frac{d-2s}{2}}}\norm{\psi}_{\dot{H}^s(Q_{N^{-1/d}}(X_i))}.
\end{equation*}
Now if we take any $\psi\in \dot{H}^s(\Tor^d)$ we have
\begin{multline*}
    \abs{\int_{\Tor^d}\psi\,d\sigma}\leq \sum_i \abs{\int_{Q_{N^{-1/d}}(X_i)}\psi\,d\sigma}\leq \sum_i \frac{1}{Nr^{\frac{d-2s}{2}}}\norm{\psi}_{\dot{H}^s(Q_{N^{-1/d}}(X_i))}\\ \les \left(\frac{1}{Nr^{d-2s}}\right)^\frac{1}{2} \left(\sum_i\norm{\psi}^2_{\dot{H}^{s}(Q_{N^{-1/d}}(X_i))}\right)^\frac{1}{2}.
\end{multline*}
To conclude, just notice that by $N\gg1$, \eqref{e:sobolevseminorm} and \eqref{e:sobolevseminorm2},
\begin{equation*}
     \sum_i\norm{\psi}^2_{\dot{H}^{s}(Q_{N^{-1/d}}(X_i))}\leq \norm{\psi}^2_{\dot{H}^{s}(\Tor^d)}.
\end{equation*}
\end{proof}
\begin{proposition}
    \label{prop:deltasinH-s} Let $d=1$ and $s>\frac{1}{2}$. Let $\mu\in\mathcal{M}^+(\Tor)$ be a sum of $N$ equispaced Dirac deltas with weight $1/N$.  Then
    \begin{equation}
        \norm{\mu-1}_{\dot{H}^{-s}}^2\lesssim \frac{1}{N^{2s}}.
    \end{equation}
\end{proposition}
\begin{proof}
    We may assume without loss of generality that
    \begin{equation*}
        \mu=\frac{1}{N}\sum_{j=0}^{N-1}\delta_{X_j},
    \end{equation*}
    where $X_{j}=j/N$. We want to compute
    \begin{equation*}
       \norm{\mu-1}_{\dot{H}^{-s}}^2= \sum_{k\in \Z\backslash\{0\}} \vert k\vert^{-2s} \vert \widehat{\mu}_k\vert^2.
    \end{equation*}
    If $k=nN$ for some $n\in \Z$, then 
    \begin{equation}
    \label{e:estmuk}
         \widehat{\mu}_k= \int_{\Tor} \exp(
    2\pi ik x)d \mu=\frac{1}{N}\sum_{j=0}^{N-1}\exp(
    2\pi i kj/N)=\frac{1}{N}\sum_{j=0}^{N-1}\exp(
    2\pi i nj)=1.
    \end{equation}
    Otherwise,   
    \begin{equation}
    \label{e:estmuk2}
        \widehat{\mu}_k= \int_{\Tor} \exp(-
    2\pi ik x)d \mu=\frac{1}{N}\sum_{j=0}^{N-1}\exp(
    2\pi i kj/N)=\frac{1-\exp(2\pi ik)}{1-\exp(2\pi i k/N)}=0,
    \end{equation}
    where in the last passage we used the equality $\sum_{i=0}^Nx^i=\frac{1-x^{N+1}}{1-x}$ for $x\neq 1$.
    Combining \eqref{e:estmuk} and \eqref{e:estmuk2} we get 
    \begin{equation*}
         \norm{\mu-1}_{\dot{H}^{-s}}^2= \sum_{k\in \Z\backslash\{0\}} \vert k\vert^{-2s} \vert \widehat{\mu}_k\vert^2=\frac{1}{N^{2s}}\sum_{n\in \Z\backslash\{0\}} \vert n\vert^{-2s}\lesssim \frac{1}{N^{2s}}.
    \end{equation*}
    
\end{proof}


\section{From local scaling laws to bounds on the lower and upper dimension of $\mu_T$}
\label{sec:dim}
In this section we recall the notion of dimension of a measure and of Ahlfors regular measures. Then, we present some results of \cite{dephilippis2023energy} that, given a minimizer $\mu \in \mathcal{A}_T^*$, link local energy bounds for its energy with the dimension of the boundary measure $\mu_T$. In \cite{dephilippis2023energy} these results were established for $s=\frac{1}{2}$ and $d=2$, here we present the results for any $d\geq 1$ and $s\in (0,1)$.
\subsection{Dimension of a measure}
\label{subsec:dimension}
Let $\alpha \in [0,d]$. If we denote by $\mathcal{H}^\alpha$ the $\alpha$-dimensional Hausdorff measure, the Hausdorff dimension of a set $A \subseteq Q$ is defined as 
\begin{equation}
    \label{def:dim_set}
    \dim A = \sup \{ \alpha \in [0, d] \colon \mathcal{H}^\alpha(A) >0\} = \inf\{ \alpha \in [0, d]\colon \mathcal{H}^\alpha(A) = 0 \}.
\end{equation}
There are many ways to define lower and upper dimension of a measure $\mu_T$ (see \cite{Mattila2000}). We mainly focus on the Hausdorff dimension, which is the notion used in \cite{dephilippis2023energy}.
\begin{definition}[Dimension of a measure] \label{def:dimension}Let $\sigma \in \mathcal{M}^+(\Tor^d).$
     The lower and upper Hausdorff dimension of $\sigma$ are defined by 
    \begin{equation}
        \label{def:hausorff_low_up}
        \underline{\textup{dim}}\,\sigma = \sup \{ \alpha \in [0, d]\colon \sigma \ll \mathcal{H}^\alpha \}  \quad \textup{ and } \quad  \overline{\textup{dim}}\,\sigma = \inf \{ \alpha \in [0, d]\colon \sigma \perp \mathcal{H}^\alpha \}.\end{equation} If the lower and upper Hausdorff dimension of $\sigma$ coincide we define the Hausdorff dimension of $\sigma$ to be the common value $\textup{dim}\,\sigma =  \underline{\textup{dim}}\,\sigma = \overline{\textup{dim}}\,\sigma . $
          
\end{definition}
Let us introduce the notion of upper $\alpha$-Ahlfors regular and lower $\alpha$-Ahlfors regular measures for $\alpha \in [0, d]$.
\begin{definition}  \label{def:more_reg_boundary}
    Let $\sigma \in \mathcal{M}^+(\Tor^d)$, and $\alpha \in [0, d].$
    \begin{itemize}
        \item[(i)] A measure is called upper $\alpha$-Ahlfors regular if there exists $M>0$ such that
        \begin{equation}
        \label{e:malphareg}
             \sigma(B_r(x)) \leq M r^\alpha \quad \textup{ for all } x \in \supp\sigma \textup{ and } r\in (0,1].
        \end{equation}
        \item[(ii)] A measure is called lower $\alpha$-Ahlfors regular if there exists $m>0$ such that
        \begin{equation}
        \label{e:malphacon}
            \sigma(B_r(x)) \geq m r^\alpha \quad \textup{ for all } x \in \supp\sigma \textup{ and } r\in (0,1].
        \end{equation} 
    \end{itemize}
    If both conditions hold (possibly for different $M,m>0$ respectively), the measure is called $\alpha$-Ahlfors regular.
\end{definition}
\begin{remark}
Notice that if $\alpha=d$ upper $\alpha$-Ahlfors regularity is equivalent to $\sigma=f dx$ for $f\in L^\infty$. Notice also that if $\alpha' < \alpha$ then 
\begin{equation*}
    \sigma \textup{ is upper } \alpha\textup{-Ahlfors regular} \Rightarrow \sigma \textup{ is upper } \alpha'\textup{-Ahlfors regular}
    \end{equation*}
    and 
    \begin{equation*}
    \sigma \textup{ is lower } \alpha'\textup{-Ahlfors regular} \Rightarrow \sigma \textup{ is lower } \alpha\textup{-Ahlfors regular}.
    \end{equation*}    
\end{remark}
By \cite[Theorem 2.5.6]{AFP} we have
\begin{proposition}
\label{prop:Malphatoabscont}
    Let $\sigma\in \M^+(\Tor^d)$. The following implications hold:
    \begin{align*}
       \sigma \textup{ is upper } \alpha\textup{-Ahlfors regular}&\Rightarrow \quad \sigma \ll \mathcal{H}^\alpha,\\ 
       \sigma \textup{ is lower } \alpha\textup{-Ahlfors regular}&\Rightarrow \quad  \sigma \perp \mathcal{H}^{\alpha'} \ \forall \alpha'>\alpha.
    \end{align*}

\end{proposition}
\begin{remark}
\label{rem:alphareganddim}
 By a direct application of the previous proposition, for upper $\alpha$-Ahlfors regular we have
    \begin{equation*}
        \underline{\textup{dim}}\,\sigma\geq \alpha, 
    \end{equation*}
    whereas for a lower $\alpha$-Ahlfors regular measure  we have
      \begin{equation*}
        \overline{\textup{dim}}\,\sigma\leq \alpha.
    \end{equation*}
    Intuitively, this means that upper $\alpha$-Ahlfors regular measures are measures that are "diffused at least $\alpha$", whereas lower $\alpha$-Ahlfors regular measures are measures that are "diffused at most $\alpha$".   

\end{remark}
Let us also recall the following characterization of the lower dimension of a measure, see \cite[Lemma 2.17]{dephilippis2023energy}.
\begin{lemma}
\label{lem:diminfcar}
    Let $\sigma\in \M^+(\Tor^d)$ with $\sigma(\Tor^d)=1$. Then
    \begin{equation*}
        \underline{\textup{dim}}\,\sigma\geq \sup \{\alpha: \alpha \leq d,  \Vert \sigma - 1\Vert_{\dot{H}^{-\frac12(d-\alpha)}} <\infty\}.
    \end{equation*}
\end{lemma}
\subsection{Local scaling laws and dimensional estimates}
Let us now recall some results from \cite{dephilippis2023energy}.
We consider a minimizer $\mu$ of $E_{s,T}$. In \cite[Theorem 1.3]{dephilippis2023energy} a first local upper bound on $I(\mu,(\T-\eps, \T))$ is proved for $d=2$. This upper bound corresponds to the worst possible scaling. It is attained when the irrigated measure is uniform, see \cite[Theorem 3.7]{dephilippis2023energy}.  The extension to arbitrary dimension is straightforward, see also the proof of Theorem \ref{thm:local_energy_lower}.
\begin{proposition}\label{prop:scaling13}
      Let $\mu \in \mathcal{A}^*_T$  be a minimizer of $E_{s,T}$ and $\eps\ll 1$. Then 
    \[ I(\mu, (\T - \eps, \T)) \lesssim  \eps^{\frac{1}{3}}.\]
\end{proposition}
We now extend the proof of   \cite[Theorem 3.5]{dephilippis2023energy} to arbitrary $d$ and $s$.
\begin{lemma}\label{thm:regularity_boundary_measure}
    Let $\mu \in \mathcal{A}_T^*$ be a symmetric minimizer of $E_{s, T}$ and assume that there exists $\beta \in (0,1)$ such that $\limsup_{\eps \rightarrow 0} \eps^{-\beta} I(\mu, (T-\eps, T)) < \infty.$ Then,
      \[ \Vert \mu_T - 1\Vert_{\dot{H}^{-\frac12(2-\alpha)}} <\infty \textup{ for all } \alpha < 2(1-s) + \frac{2\beta}{1+\beta}.\]
\end{lemma}
\begin{proof} Let $\eta \in (0,1)$, $\eps\in (0,T)$ and consider $\mu_t = \sum_{i} \varphi_i \delta_{X_i(t)}$. Let $e_j$ for $j=1,\dots, d$, denote the canonical base of $\R^d$. Let us build a competitor for $\mu$. Since $\mu$ is symmetric we can restrict to $(0,T)$. Let
    \[X_i^\pm(t) =  \begin{cases} 
    X_i(t) \quad & \textup{ if } t\in (0, T-\eps),\\
    X_i(t) \pm \frac\eta\eps (t- (T-\eps))e_j \quad & \textup{ if }t\in ( T-\eps,T).
    \end{cases}\] 
     We set $\mu^\pm_t=\sum_{i} \varphi_i \delta_{X_i^\pm(t)}$ and define $\tilde{\mu}= \frac12(\mu^+ + \mu^-)$. We have by a direct computation 
    \begin{equation}
        \label{prop:energy_competitor}
        I(\tilde{\mu}, (T-\eps,T)) \lesssim I(\mu, (T-\eps,T)) + \frac{\eta^2}{\eps}. 
    \end{equation}
    By minimality we therefore conclude
    \begin{eqnarray}
        \label{est:difference-boundary-norms}
        && \Vert \mu_T -1 \Vert_{\dot{H}^{-s}}^2 - \Vert \tilde{\mu}_T -1 \Vert_{\dot{H}^{-s}}^2\leq  I(\tilde{\mu}, (T-\eps,T)) - I(\mu, (T-\eps,T))\nonumber \\
        && \lesssim I(\mu,(T-\eps, T)) + \frac{\eta^2}{\eps}\lesssim_\mu \eps^\beta + \frac{\eta^2}{\eps} \lesssim_\mu \eta^{\frac{2\beta}{1+\beta}} 
    \end{eqnarray}
    choosing $\eps = \eta^{2/(1+\beta)}$ to minimize the right-hand side of \eqref{est:difference-boundary-norms}. Now, we turn our attention to the left-hand side of \eqref{est:difference-boundary-norms}. We set $\sigma = \mu_T -1$ and $\tilde{\sigma}= \tilde{\mu}_T -1$. Notice that if we set $\tau_{\eta}(x)=x-\eta e_j$ we have
    \begin{equation*}
        \Tilde{\sigma}=\frac{1}{2}(\tau_\eta\#\sigma+\tau_{-\eta}\#\sigma).
    \end{equation*}
    This implies 
    \[ \widehat{\tilde{\sigma}}_k = \int_{\Tor^d} \exp(i2\pi k\cdot x) \,d\tilde{\sigma}(x)  = \cos(2\pi \eta k\cdot e_j) \widehat{\sigma}_k.\]
    This leads us to 
    \begin{equation}
        \label{est:differenc-norms:below}
         \Vert \mu_T -1 \Vert_{\dot{H}^{-s}}^2 - \Vert \tilde{\mu}_T -1 \Vert_{\dot{H}^{-s}}^2 = \sum\limits_{k\in \Z^d} \vert k\vert^{-2s} \sin^2(2\pi k\cdot e_j\eta)\vert \widehat{\sigma}_k \vert^2 \gtrsim \eta^2 \sum_{\vert k \vert \lesssim \eta^{-1}} \vert k\cdot e_j\vert^2 \vert k\vert^{-2s}\vert \widehat{\sigma}_k\vert^2 .
     \end{equation}
     Since this argument holds for any $j=1, \dots, d$, it follows from \eqref{est:difference-boundary-norms} that
     \begin{equation*}
         \eta^{\frac{2\beta}{1+\beta}}\gtrsim_\mu \eta^2 \sum_{\vert k \vert \lesssim \eta^{-1}}  \vert k\vert^{2(1-s)}\vert \widehat{\sigma}_k\vert^2.
     \end{equation*}
Using the characterization \eqref{e:H-scharacterization1} and this last estimate yields for $\alpha >0$,
\begin{eqnarray}
    && \Vert \sigma \Vert_{\dot{H}^{-\frac12(2-\alpha)}}^2  \lesssim_\mu \int_0^1 \eta^{2\frac12 (2-\alpha)+ 2(1-s)}\left( \sum_k \vert k \vert^{2(1-s)} \vert \widehat{\sigma}(k)\vert^2 \right)\frac{\,d\eta}{\eta }\nonumber \\
    && \lesssim_\mu \int_0^1 \eta^{ 2(1-s)-\alpha } \eta^{\frac{2\beta}{1+\beta}}\frac{\,d\eta}{\eta} .\label{proof:regularity_estimate}
\end{eqnarray}
The integral on the right-hand side of \eqref{proof:regularity_estimate} is finite if and only if $ 2(1-s)-\alpha + \frac{2\beta}{1+\beta}> 0$. Hence, the claim follows.
\end{proof}

 We can now apply the strategy from \cite[Corollary 3.6 \& Theorem 3.7]{dephilippis2023energy} to obtain the following relation between local scaling laws and the dimension of $\mu_T$.
\begin{proposition}
\label{prop:diminfbound}
    Let $s\in (0,1)$ and $\mu\in \mathcal{A}^*_{T}$ be a symmetric minimizer of $E_{s, T}$.  Assume that for some $\beta \in (0,1)$ with $s>\beta/(1+\beta)$\begin{equation*}\limsup_{\eps \rightarrow 0} \eps^{-\beta} I(\mu, (\T-\eps, \T)) < \infty.\end{equation*} Then
    \begin{equation}
        \label{e:dimbound}
        d-2s+\frac{2\beta}{1+\beta}\leq \underline{\textup{dim}}\,\mu_T\leq \overline{\textup{dim}}\,\mu_T \leq \frac{2d(1-\beta)}{1 + \beta }.
    \end{equation}
\end{proposition}
\begin{proof}
    The proof is basically the same as in \cite{dephilippis2023energy}, in which the result is proved for $s=\frac{1}{2}$ and $d=2$. In particular, the first estimate is a direct consequence of  Lemma \ref{lem:diminfcar} combined with Lemma \ref{thm:regularity_boundary_measure}. The second estimate does not depend on $s$ and the proof found in \cite[Theorem 3.7]{dephilippis2023energy} generalizes directly to any dimension.
\end{proof}
\begin{remark}
\label{rem:dimension}
   In the last proposition we assume that 
   \begin{equation*}
       s>\frac{\beta}{1+\beta}.
   \end{equation*}
   If this does not hold, Lemma \ref{thm:regularity_boundary_measure} together with  Lemma \ref{lem:diminfcar} directly imply that
   \begin{equation*}
       \underline{\textup{dim}}\,\mu_T=\overline{\textup{dim}}\,\mu_T=\textup{dim}\,\mu_T=d.
   \end{equation*}
   In particular, since by Proposition \ref{prop:scaling13}, we have  $\beta \ge 1/3$, we may restrict ourselves to $s>1/4$.
\end{remark}

\begin{remark}\label{rem:critical}
    Equality of both sides in \eqref{e:dimbound} is achieved for
    \begin{equation}
        \label{e:critical}
         \beta_{c}(s, d) = \frac{d+2s}{3d+2(1-s)}.
    \end{equation}
    This would correspond to
    \[ \textup{dim}\,\mu_T = \frac{2d(d+1-2s)}{2d+1}.\]
     The critical exponent \eqref{e:critical} is equal to the exponent of the global scaling law for $0<T\leq 1$ and $d>2s$, see Proposition \ref{pro:global_scaling}. In particular, for $s=1/2$ and $d=2$  we get $\beta_{c}(1/2,2) = 3/7$, corresponding to the conjectured $\textup{dim}\,\mu_T = 8/5$. For $d=1$ both sides of \eqref{e:dimbound} coincide for $\beta_c(s,1) =(1+2s)/(5-2s)$ and correspondingly, $\textup{dim}\,\mu_T= 4/3(1-s)$.
\end{remark}
In the subsequent sections we will prove Theorem \ref{thm:lower_local_bound} and Theorem \ref{thm:local_energy_lower}, which relates upper $\alpha$-Ahlfors regularity and lower $\alpha$-Ahlfors regularity of $\mu_T$ with local scaling laws. We postpone the proofs of these results and first show  how they combine with Proposition \ref{prop:diminfbound} to  prove Theorem \ref{thm:Malphareg}.
\begin{proof}[Proof of Theorem \ref{thm:Malphareg}]
Assume that $\mu_T$ is upper $\alpha$-Ahlfors regular. We can suppose $s>\frac{1-\alpha}{2}$ otherwise there is nothing to prove. Then, by Theorem \ref{thm:lower_local_bound},  for every $\delta\in(0,\beta_{\textup{reg}})$ (recall the definition of $\beta_\textup{reg}$ \eqref{e:betareg})
\begin{equation*}
    I(\mu,(T-\eps,T))\lesssim_\delta \eps^{\beta_{\textup{reg}}(s,\alpha)-\delta}.
\end{equation*}
 By Proposition \ref{prop:diminfbound} and Remark \ref{rem:alphareganddim} we have
\begin{equation*}
      \alpha\leq \underline{\dim}\,\mu_T\leq \overline{\dim}\,\mu_T\leq \frac{2(1-\beta_{\textup{reg}}(s,\alpha)+\delta)}{1 +\beta_{\textup{reg}}(s,\alpha)-\delta }.
\end{equation*}
This yields recalling \eqref{e:alphabar}
\begin{equation*}
    \alpha\leq \Bar{\alpha}+C(s,\alpha)\delta
\end{equation*}
    for some $C(s,\alpha)>0$ and all $\delta\in(0,\beta_{\textup{reg}})$.  Letting $\delta\rightarrow 0$, this gives \eqref{e:alphaless} .
    Analogously, assume that $\mu_T$ is lower $\alpha$-Ahlfors regular. Then, by Theorem \ref{thm:local_energy_lower} we have (recall the definition \eqref{e:betacon} of $\beta_{\textup{con}}(\alpha)$)
    \begin{equation*}
         I(\mu, (T-\eps, T)) \lesssim_\mu \eps^{\beta_{\textup{con}}(\alpha)}.
    \end{equation*}
    Thus, by Proposition \ref{prop:diminfbound} we have 
    \begin{equation*}
        1-2s+\frac{2\beta_{\textup{con}}(\alpha)}{1+\beta_{\textup{con}}(\alpha)}\leq \underline{\textup{dim}}\,\mu_T\leq \overline{\textup{dim}}\,\mu_T \leq \alpha
    \end{equation*}
    which yields $\alpha\geq \Bar{\alpha}$ using the definitions \eqref{e:alphabar} and \eqref{e:betacon} of $\Bar{\alpha}$ and  $\beta_{\textup{con}}$ respectively.
\end{proof}
    






\section{Local scaling laws for minimizers with upper $\alpha$-Ahlfors regular boundary measures}
\label{sec:regular}
The goal of this section is to prove Theorem \ref{thm:lower_local_bound}. This section is divided in two parts. In the first part (which holds for any $d\ge 1$), we prove that the class of   upper $\alpha$-Ahlfors regular measures  is stable  with respect to  McCann interpolation. This result could be of independent interest.  We then show that upper $\alpha$-Ahlfors regular measures are embedded in some  Sobolev space of negative order.  In the second part we fix $d=1$ and consider a minimizer $\mu$ of $E_{s,T}$. We prove that if $\mu_T$ is upper $\alpha$-Ahlfors regular, then $\mu$ satisfies a certain local scaling law.
\subsection{Interpolated measures, upper $\alpha$-Ahlfors regularity and Sobolev regularity}
\label{subsec:interpolatedsobolev}
\begin{lemma}[Upper $\alpha$-Ahlfors regularity of displacement measures]
    \label{lem:reg_displacement} Let $\pi\in \M^+_{{\rm loc}}(\R^d\times\R^d)$ be a monotone transport plan i.e. for every $(x_1,y_1),(x_2,y_2)\in \supp\pi$
    \begin{equation}
        \label{e:monotonicity}
        (x_1-x_2)\cdot (y_1-y_2)\geq 0
    \end{equation}
     and for $\lambda\in (0,1)$ let $\mu^\lambda=F_\lambda\# \pi$ where  $F_\lambda(x,y)=(1-\lambda) x+\lambda y$. If  $\mu^0=\pi_1$ is upper $\alpha$-Ahlfors regular for some $M>0$ then $\mu^\lambda$ is upper $\alpha$-Ahlfors regular for some constant $M_{\lambda}\lesssim (1-\lambda)^{-\alpha} M$.
\end{lemma}
\begin{proof}
    We claim that for every $r\in(0,1)$ and $x_{0}\in \R^d$ there exists $x_{\lambda,r}$ such that
    \begin{equation}
        F_\lambda^{-1}(B_{r}(x_{0}))\cap \supp{\pi}\subset B_{2r/(1-\lambda)}(x_{\lambda,r})\times \R^d. \label{subset:reg_disp}
    \end{equation} 
    Indeed, this would imply
    \begin{equation*}
        \mu^{\lambda}(B_{r}(x_{0}))=\pi(F_\lambda^{-1}(B_{r}(x_{0})))\leq \pi(B_{2r/(1-\lambda)}(x_{\lambda,r})\times \R^d)=\mu^0(B_{2r/(1-\lambda)}(x_{\lambda,r}))\leq M2^\alpha(1-\lambda)^{-\alpha}r^\alpha,
    \end{equation*}
    which would conclude the proof of the upper $\alpha-$Ahlfors regularity of $\mu^\lambda$.
    We now prove claim \eqref{subset:reg_disp}. To this aim, it is enough  to show that for any $(x_1,y_1),(x_2,y_2)\in F_\lambda^{-1}(B_{r}(x_{0}))\cap \supp{\pi}$ we have $|x_1-x_2|\leq 2r/(1-\lambda)$. By definition, for any $(x_1,y_1),(x_2,y_2)\in F_\lambda^{-1}(B_{r}(x_{0}))$ we have
    \begin{equation*}
     |(1-\lambda)(x_1-x_2)+\lambda(y_1-y_2)|=|F_\lambda(x_1,y_1)-F_\lambda(x_2,y_2)|\leq 2r.
    \end{equation*}
   By \eqref{e:monotonicity}, we then have
    \begin{equation}
    \label{e:usemonotonicity}
        2r\geq |(1-\lambda)(x_1-x_2)+\lambda(y_1-y_2)|\geq (1-\lambda) |x_1-x_2|
    \end{equation}
    and the claim is proved.
\end{proof}
\begin{remark}
Notice that the measure $\mu^\lambda$ can be in general much less concentrated than $\mu^0$ (if for instance $\mu^1=\pi_2$ is upper $\alpha'$-Ahlfors regular for some  $\alpha'>\alpha$).
\end{remark}
\begin{remark}
    The arguments in the proof of Lemma \ref{lem:reg_displacement} are reminiscent of the ones of \cite[Theorem 8.7]{Villani} where it is shown that absolute continuity with respect to the Lebesgue measure is preserved by displacement interpolation. Lemma \ref{lem:reg_displacement} can also be seen as a generalization of  the $L^\infty$ bounds on the McCann interpolant \eqref{e:mccann}, see \cite[Corollary 19.5]{Villani}. Our proof is however  quite elementary compared to the standard proof based on the Monge-Ampère equation and the log-concavity of the determinant. Let us however point out that as opposed to \cite[Theorem 19.4]{Villani}, we do not get sharp constants (which are essential in many geometric applications of \cite[Theorem 19.4]{Villani}).
\end{remark}

\begin{proposition}[Sobolev embedding] 
    \label{prop:Mksobolevemb}
    Assume that $\mu\in \mathcal{P}(\Tor^d)$ is upper $\alpha$-Ahlfors regular for some $M>0$. Then $\mu-1\in \dot{H}^{-\gamma}(\Tor^d)$ for every $\gamma$ such that 
    \begin{equation}
    \label{condWgammap}
    \gamma>\frac{d-\alpha}{2}.
    \end{equation}
    Moreover,
        \begin{equation*}
        \|\mu-1\|^2_{\dot{H}^{-\gamma}}\lesssim_\gamma M+1.
    \end{equation*}
\end{proposition}
\begin{proof}
    Consider a family of convolution kernels $(\rho_\eps)_\eps,\, \eps>0$. We want to use \eqref{e:H-scharacterization2}. For every $x\in\Tor^d$ we estimate 
    \begin{equation}
        \label{e:convest}  
    \rho_\eps\ast \mu(x)=\int_{B_\eps(x)} \rho_\eps(x-y)d\mu(y)
    \lesssim \frac{\mu(B_\eps(x))}{\eps^d}.
    \end{equation}
    Moreover,  
    \begin{equation}
        \label{e:estmu}
\int_{\Tor^d} \mu(B_\eps(x))^2\,dx\lesssim \sup_{x\in\Tor^d} \mu(B_\eps(x))\int_{\Tor^d} \mu(B_\eps(x))\,dx\lesssim M \eps^{\alpha} \eps^d,
    \end{equation}
    where in the last passage we used upper $\alpha$-Ahlfors regularity and Fubini's Theorem.
   Notice finally that 
    \begin{equation}
    \label{e:plusone}
    \abs{\rho_\eps\ast(\mu-1)(x)}^2\lesssim \abs{\rho_\eps\ast\mu(x)}^2+1.
    \end{equation}
    Plugging \eqref{e:convest}, \eqref{e:estmu} and \eqref{e:plusone} in \eqref{e:H-scharacterization2} we obtain
    \[
    \|\mu-1\|_{\dot{H}^{-\gamma}}^2\lesssim  (M+1)\int_0^1 \frac{\eps^{2\gamma}}{\eps^{2d}}\eps^{\alpha} \eps^d \frac{d\eps}{\eps} = (M+1)\int_0^1 \eps^{2\gamma + \alpha - d -1} d \eps.
    \]
    This integral is finite if and only if 
    \[
    2\gamma +\alpha > d,
    \]
    which yields the condition \eqref{condWgammap}.
\end{proof}

\subsection{Local scaling laws for minimizers with upper $\alpha$-Ahlfors regular boundary measure }
In this section we fix $s>\frac{1}{4}$.  Let $\mu$ be a minimizer of $E_{s,T}$. Goal of this section is to prove Theorem \ref{thm:lower_local_bound}. We start with the following lemma.
\begin{lemma}
    \label{lem:H-s}
    Let $\mu,\nu\in \mathcal{P}(\Tor)$ be such that $\mu$ is upper $\alpha$-Ahlfors regular for some $M>0$ and  $\alpha>d-2s$. Let $\pi$ be an optimal transport plan for $W(\mu,\nu)$. For $\lambda\in (0,1)$ let $\mu^\lambda=F_\lambda\# \pi$ where  $F_\lambda(x,y)=(1-\lambda) x+\lambda y$. Then, letting
    \begin{equation}
    \Bar{\theta}(s,\alpha)=s-\frac{d-\alpha}{2},
\end{equation}
we have for every $\theta<\Bar{\theta}(s,\alpha)$,
\begin{equation}
\label{e:hsinterpestimate}
     \Vert \mu- \mu^\lambda \Vert_{\Hdots} \lesssim_{\lambda,\theta,\mu} W^\theta(\mu, \nu).
\end{equation}
\end{lemma}
\begin{proof}
     Fix $\eps\in(0,\eps_0)$ and consider a regular convolution kernel $\rho_\eps$, recall \eqref{prop:reg_conv_kernel}. Using triangle inequality we have    
    \begin{equation}
    \label{e:triangleineq}
        \Vert \mu - \mu^\lambda\Vert_{\Hdots} 
         \leq \Vert \mu - \mu \ast \rho_\eps \Vert_{\Hdots} + \Vert \mu \ast \rho_\eps - \mu^\lambda\ast \rho_\eps \Vert_{\Hdots} + \Vert \mu^\lambda - \mu^\lambda \ast \rho_\eps \Vert_{\Hdots}.    
    \end{equation}

Let us estimate the first term on the right-hand side of the previous inequality. By Proposition \ref{prop:Mksobolevemb} we have $ \|\mu-1\|_{\dot{H}^{-\gamma}}^2\lesssim_\gamma M+1$  for any $\gamma>\frac{d-\alpha}{2}$. Moreover by Proposition \ref{prop:kernelestimate}
we have  
\begin{equation*}
    \Vert \mu- \mu \ast \rho_\eps \Vert_{\Hdots}\lesssim \eps^{s-\gamma} \Vert \mu-1 \Vert_{\dot{H}^{-\gamma}}\lesssim_{\gamma,\mu} \eps^{s-\gamma} (M+1)^{\frac{1}{2}}.
\end{equation*}
We now turn to the third right-hand side term in \eqref{e:triangleineq}. Since $\pi$ is monotone, by Lemma \ref{lem:reg_displacement}, $\mu^\lambda$ is upper $\alpha$-Ahlfors regular for some $M_{\lambda}\lesssim_\lambda M$, so that arguing exactly as above,
\begin{equation*}
    \Vert \mu^\lambda- \mu^\lambda \ast \rho_\eps \Vert_{\Hdots}\lesssim_\lambda \eps^{s-\gamma} \Vert \mu-1 \Vert_{\dot{H}^{-\gamma}}\lesssim_{\lambda,\gamma,\mu} \eps^{s-\gamma} (M+1)^{\frac{1}{2}}.
\end{equation*} 
Let us now focus on the second right-hand side term in \eqref{e:triangleineq}. Using H\"older inequality with $p=1/(1-s)$ and $q=1/s$ we find that for any $\sigma\in \dot{H}^{-s}(\Tor)$
\begin{equation*}
\Vert \sigma \Vert_{\Hdots}=\left(\sum_{k\in\Z}\abs{k}^{-2s}\abs{\widehat{\sigma}_k}^{2}\right)^\frac{1}{2}=\left(\sum_{k\in\Z}\abs{\widehat{\sigma}_k}^{2(1-s)}\abs{k}^{-2s}\abs{\widehat{\sigma}_k}^{2s}\right)^\frac{1}{2} \leq \Vert \sigma \Vert_{L^2}^{1-s}\Vert \sigma \Vert_{\dot{H}^{-1}}^s.
\end{equation*}
Hence,

\begin{equation*}
\begin{split}
    & \Vert \mu \ast \rho_\eps - \mu^\lambda\ast \rho_\eps \Vert_{\Hdots} \leq \Vert \mu \ast \rho_\eps - \mu^\lambda\ast \rho_\eps \Vert_{L^2}^{1-s}\Vert \mu \ast \rho_\eps - \mu^\lambda\ast \rho_\eps \Vert_{\dot{H}^{-1}}^s\\
    &\quad  \leq (\Vert \mu \ast \rho_\eps \Vert_{L^2} + \Vert \mu^\lambda \ast \rho_\eps \Vert_{L^2})^{1-s} \Vert \mu \ast \rho_\eps - \mu^\lambda\ast \rho_\eps \Vert_{\dot{H}^{-1}}^s. 
     \end{split}
\end{equation*}
Moreover, by Proposition \ref{prop:H-sW2},
\begin{equation*}
    \label{e:H-1estimate}
    \Vert \mu \ast \rho_\eps - \mu^\lambda\ast \rho_\eps \Vert_{\dot{H}^{-1}}\lesssim \max\{\Vert \mu \ast \rho_\eps \Vert_{L^\infty}^{\frac{1}{2}},\Vert \mu^\lambda \ast \rho_\eps \Vert_{L^\infty}^{\frac{1}{2}}\} W_{\textup{per}}(\mu\ast \rho_\eps, \mu^\lambda \ast \rho_\eps).
\end{equation*}
Now, arguing as in the proof of Proposition \ref{prop:Mksobolevemb} we get 
\begin{equation*}
    \rho_\eps\ast \mu(x)
    \les \frac{\mu(B_\eps(x))}{\eps^d}\le \frac{M}{\eps^{d-\alpha}}.
\end{equation*}
Hence 
\begin{equation*}
    \Vert \mu \ast \rho_\eps \Vert_{L^2}^2\leq \Vert \mu \ast \rho_\eps \Vert_{L^\infty} \int_{\Tor}(\mu \ast \rho_\eps)(x)\,dx \lesssim  \frac{M}{\eps^{d-\alpha}}.
\end{equation*}
Moreover, again by Lemma \ref{lem:reg_displacement} the same estimates hold also for $\mu^\lambda$, that is
\begin{equation*}
        \Vert \mu^\lambda \ast \rho_\eps \Vert_{L^2}^2 \leq \Vert \mu^\lambda \ast \rho_\eps \Vert_{L^\infty}\lesssim_\lambda  \frac{M}{\eps^{d-\alpha}}.
\end{equation*}
Finally, by \cite[Lemma 5.2]{santambrogio2015optimal},
\begin{equation*}
    W_{\textup{per}}(\mu\ast \rho_\eps, \mu^\lambda \ast \rho_\eps)\leq W_{\textup{per}}(\mu, \mu^\lambda)\le W(\mu,\mu^\lambda)\stackrel{\eqref{e:geodesic}}{=}\lambda W(\mu, \nu).
\end{equation*} 
Putting all the estimates together we get
\begin{equation*}
     \Vert \mu - \mu^\lambda
     \Vert_{\Hdots} \lesssim_{\lambda,\gamma,\mu} \eps^{s-\gamma}+\eps^{\frac{\alpha-d}{2}}W^s(\mu, \nu).
\end{equation*}
We optimize in $\eps$ choosing 
\begin{equation*}
    \eps= W^{\frac{s}{s-\beta}}(\mu, \nu) \quad \textup{ where } \quad \beta = \gamma-\frac{d-\alpha}{2}
\end{equation*}
if $W^{\frac{s}{s-\beta}}(\mu, \nu)\leq \eps_0$ and $\eps=\eps_0$ otherwise.
In both cases this implies 
\begin{equation*}
    \Vert \mu - \mu^\lambda
     \Vert_{\Hdots} \lesssim_{\lambda,\theta,\mu} W^{s-\frac{s(d-\alpha)}{2(s-\beta)}}(\mu, \nu).
\end{equation*}
Since this holds for any $\beta>0$, we conclude the proof of \eqref{e:hsinterpestimate}.
\end{proof}
\begin{proof}[Proof of Theorem \ref{thm:lower_local_bound}]
 Let $\pi$ be an optimal transport plan for $W(\mu_{T-\eps},\mu_T)$. \\

    \quad \emph{\textbf{Step 1: Construction of a competitor. }}
     Let $X$ be the Lagrangian representation of $\mu$ given by Proposition \ref{prop:lagrangian}. We first notice that since $d=1$ by Remark \ref{rem:periodicnonperiodic} and Lemma \ref{lem:monotonecoupling}, the map $\psi=X(\cdot,T-\eps)$ is monotone and thus $\pi=({\rm Id}\times \psi)\#\mu_T$ as well. Since $d=1$, it is thus also an optimal transport plan for $W^2(\mu_{T-\eps},\mu_T)$. We then set
     \begin{equation*}
         Y(x,t)=
         \begin{cases}
             X(x,t) & \textup{ for } t\leq T-\eps \\
             \frac{1}{2}(X(x,T-\eps)+X(x,t)) & \textup{ for } T-\eps < t\leq T,
         \end{cases}
     \end{equation*}
     and $\nu_t=Y(\cdot,t)\# \mu_T$. With the notation of Lemma \ref{lem:H-s}, $\nu_T=F_{1/2}\#\pi=\mu_T^{1/2}$. Notice that  if $\mu_t= \sum_i \varphi_i \delta_{X_i(t)}$, this construction is equivalent to setting $\nu_t=\sum_{i} \varphi_i \delta_{Y_i(t)}$ where
 \begin{equation*}
 Y_{i}(t)=
     \begin{cases}
         X_{i}(t) &\textup{ if } t\leq T-\eps \\
          \frac{1}{2}(X_i(t)+X_i(T-\eps)) &\textup{ if } t\leq T-\eps.
     \end{cases}
     \end{equation*}

          \begin{figure}[h]
         \centering
         \includegraphics[scale = 0.6]{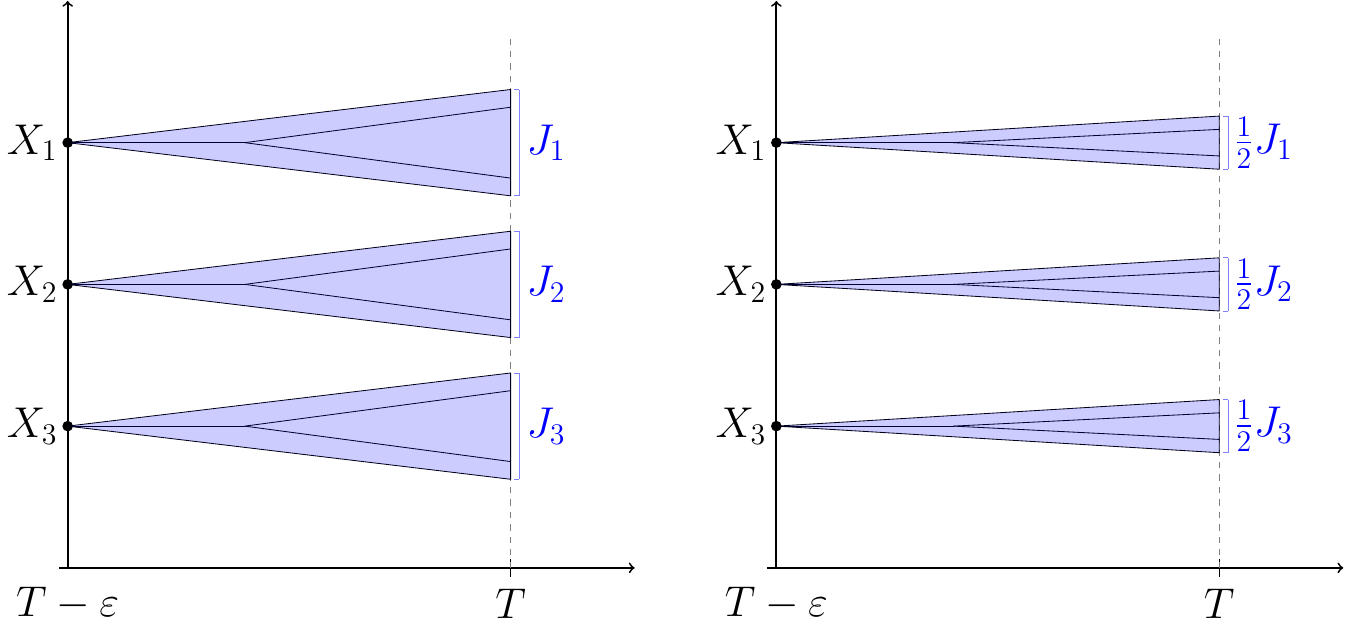}
         \caption{On the left a sketch of $\mu$ on $[T-\eps,T]$, on the right a sketch of $\nu$ on $[T-\eps,T]$. }
         \label{fig:slowing_down_construction}
     \end{figure}

         \quad \emph{\textbf{Step 2: Estimate of the internal energy of $\mu_T$.} }\\  Let us compare  the minimizer $\mu$ to the competitor we constructed in Step 1. By construction  we have (notice that for $t\ge T-\eps$, $\Dot{Y}_{i}(t)= \Dot{X}_{i}(t)/2$)
         \begin{equation*}
         \begin{split}
             &I(\nu, (0, T-\eps))=I(\mu, (0, T-\eps)), \quad  P(\nu, (T-\eps, T))= P(\mu, (T-\eps, T)) \\
            & \quad \textup{ and } \quad  
          E_{\textup{kin}}(\nu, (T-\eps, T))= \frac{1}{4}E_{\textup{kin}}(\mu, (T-\eps, T)).
          \end{split}
         \end{equation*}
         The minimality of $\mu$ yields
  \begin{equation}
      \label{Est:internal_energy}
      \frac{3}{4} E_{\textup{kin}}(\mu, (T-\eps, T)) \leq   \Vert \nu_T - 1\Vert_{\Hdots}^2 - \Vert \mu_T -1 \Vert_{\Hdots}^2 = \Vert \mu_T - \nu_T \Vert^2_{\Hdots} + 2\langle \mu_T - \nu_T, \mu_T - 1\rangle_{\Hdots}.
  \end{equation}
     Using Lemma \ref{lem:H-s} we get 
\begin{equation*}
     \Vert \mu_T - \nu_T \Vert_{\Hdots} \lesssim_{\theta,\mu} W^{\theta}(\mu_T, \mu_{T-\eps}) \quad \textup{ for } \theta< s-\frac{1-\alpha}{2}.
\end{equation*}
Let us now focus on the scalar product in \eqref{Est:internal_energy}. By Proposition \ref{prop:Mksobolevemb} we have $\mu_T-1\in \dot{H}^{-\gamma}(\Tor)$ for all $\gamma>\frac{1-\alpha}{2}$ with $\Vert \mu_T-1\Vert_{\dot{H}^{-\gamma}}^2\les_{\gamma} M+1$.  We then have by H\"older inequality
\begin{multline*}
    \langle\mu_T-\nu_T,\mu_T-1\rangle_{\dot{H}^{-s}}=\sum_{k\in \mathbb{Z}}|k|^{-2s}\widehat{(\mu_T-\nu_T)}_k\overline{\widehat{(\mu_T-1)}}_k= \sum_{k\in \mathbb{Z}}|k|^{-2s+\gamma}\widehat{(\mu_T-\nu_T)}_k \cdot |k|^{-\gamma}\overline{\widehat{(\mu_T-1)}}_k\\ 
    \leq \left(\sum_{k\in \mathbb{Z}}|k|^{-4s + 2\gamma}|\widehat{\mu_T-\nu_T}|^2_k\right)^\frac{1}{2}\left( \sum_{k\in \mathbb{Z}}|k|^{-2\gamma}|\widehat{\mu_T-1}|^2_k\right)^\frac{1}{2}=\Vert \mu_T-\nu_T\Vert_{\dot{H}^{-(2s -\gamma)}} \Vert \mu_T-1\Vert_{\dot{H}^{-\gamma}}.
\end{multline*}
From \eqref{e:hsinterpestimate} again we obtain 
\begin{equation*}
\Vert \mu_T-\nu_T\Vert_{\dot{H}^{-(2s -\gamma)}} \les_{\tilde{\theta},\mu} W^{\tilde{\theta}}(\mu_{T-\eps}, \mu_T) \quad \textup{ for } \tilde{\theta} < 2s -\gamma -\frac{1-\alpha}{2}.
\end{equation*}
Now fix $\delta>0$ and choose 

\begin{align*}
   &\gamma=\frac{1-\alpha}{2}+\delta,\\
    &\theta= s-\frac{1-\alpha}{2}-\frac{\delta}{2},\\
    &\Tilde{\theta}=2s-\gamma-\frac{1-\alpha}{2}.
\end{align*}

Putting all the estimates together we arrive at 
\begin{multline}
\label{e:scalingdeltaxi}
     E_{\textup{kin}}(\mu, (T-\eps, T)) \lesssim_{\delta,\mu} (W^{s-\frac{1-\alpha}{2}-\frac{\delta}{2}}(\mu_T, \mu_{T-\eps}))^2 + \Vert \mu_T-1\Vert_{\dot{H}^{-\gamma}}W^{2s -\gamma -\frac{1-\alpha}{2}-\delta}(\mu_T, \mu_{T-\eps}) \\ \lesssim_{\delta,\mu}  W^{2s-(1-\alpha)-\delta}(\mu_T, \mu_{T-\eps}).
\end{multline}
Since $E_{\textup{kin}}(\mu, (T-\eps, T)) \geq \frac{1}{\eps}W^2(\mu_T, \mu_{T-\eps})$ by \eqref{e:ecinestimate} (see Remark \ref{rem:periodicnonperiodic}) and since \eqref{e:scalingdeltaxi} holds for any $\delta>0$, we obtain after simplification that for every $\beta<\beta_{\textup{reg}}$ (recall \eqref{e:betareg}),
\begin{equation*}
    E_{\textup{kin}}(\mu,(T-\eps,T))\lesssim_{\beta,\mu}  \eps^\beta.
\end{equation*}
 Finally, by
Proposition \ref{prop:equipartition} integrated on $[T-\eps,T]$ we have
\begin{equation*}
   P(\mu,(\T-\eps,\T)) =E_{\textup{kin}}(\mu,(\T-\eps,\T))+\Lambda\eps.
\end{equation*}
This implies
\begin{equation*}
    P(\mu,(\T-\eps,\T))\les_{\beta,\mu} \eps^{\beta}+\Lambda\eps \lesssim_{\beta,\mu} \eps^{\beta},
\end{equation*}
thus the estimate holds for $I(\mu, (T-\eps, T))$ as well. This concludes the proof of Theorem \ref{thm:lower_local_bound}.
\end{proof}
\begin{remark}
    Most of the arguments of the proof also hold  for $d>1$. The hypothesis $d=1$ plays a role only to ensure that the coupling $\pi$ between $\mu_{T-\eps}$ and $\mu_T$ induced by a minimizer is the optimal transport plan for $W$, so that Lemma \ref{lem:H-s} can be applied.
\end{remark}

\section{Local scaling laws for minimizers with lower $\alpha$-Ahlfors regular boundary measures}
\label{sec:concentrated}
The aim of this section is to prove Theorem \ref{thm:local_energy_lower}.
We start by adapting the construction from \cite{dephilippis2023energy} to measures with small support. This is then used to construct a competitor that irrigates with small energy any boundary measure $\mu_T$ which is lower $\alpha$-Ahlfors regular (recall Definition \ref{def:more_reg_boundary}). Let us insist on the fact that as opposed to the previous section, we do not modify  the boundary condition $\mu_T$ for $t=T$.

\subsection{Local construction for concentrated boundary measures}
In this section, we provide a local construction for boundary measures with small support. We make the construction for arbitrary $d\ge 1$. 
\begin{lemma}\label{lem:single_branch}
Let $\eps, \Phi>0$, $r\in (0,1)$ and  $\mu^{\pm}\in \mathcal{M}^+(\R^d)$ be such that  $\mu^\pm(\R^d)=\Phi$. Assume that there exist $x_\pm \in \R^d$ such that  $ \supp(\mu^\pm)\subseteq B_{r}(x_{\pm})$. Then, there is $\mu \in \mathcal{A}^*(\R^d,(0,\eps))$ and $C>0$ such that 
\begin{equation}
\label{e:localconstruction}
I(\mu, (0,\eps)) \leq \frac{1}{\eps} W^2(\mu^-, \mu^+) + C\left(\frac{r^2}{\eps}\Phi + \eps\Phi^\frac{d-1}{d}\right) \end{equation}
and
\begin{equation}\label{e:supmut}\supp(\mu_t) \subseteq \frac{\eps - t}{\eps} B_r(x_{-}) + \frac{t}{\eps} B_r(x_{+}) \quad \textup{ for all } t \in [0, \eps].\end{equation}
\end{lemma}
\begin{proof} The construction follows the lines of  \cite[Proposition 3.9]{dephilippis2023energy}. The main difference is that in our case we have a bound on the  diameter of the support of the measures.\\

 \quad \textbf{\emph{Step 1: Construction on $[0,2r]^d$}}
Suppose first $x_-=x_+$. Up to translation, we may assume that $x_\pm=0$ and thus  $\supp{\mu^\pm}\subseteq [0,2r]^d$. Since $\mu^+(\R^d) = \mu^-(\R^d)$, consider an optimal transport plan $\overline{\pi}$ for $W(\mu^-,\mu^+)$ and set
\begin{equation}
    \label{building_block:linear_interpolation}
    \int_{\R^d} \psi \,d\nu^t= \int_{\R^d\times \R^d} \psi\left(\frac{\eps-t}{\eps} x + \frac{t}{\eps} y\right) \,d\overline{\pi}(x,y) \quad \textup{ for } \psi \in C^0(\R^d) \textup{ and } t\in [0, \eps].
\end{equation}
We only perform the construction on $[0, \eps/2]$. The construction on $[\eps/2,\eps]$ is done analogously.
Let $\delta \in (1/4, 1/2)$ and  for every $k\geq 0 $ let $N_k=2^k$,
\begin{equation}
    \label{discrete_time_steps_construction}
     t_k =  \delta^k\frac\eps2 \qquad \textrm{so that } \qquad  t_k - t_{k+1} = \delta^{k}(1-\delta) \frac\eps2.
\end{equation}
 We divide $[0,2r]^d$ into $N_k^d$ cubes $Q_{i,k}$ of side length $l_k=2rN_k^{-1}$. We denote by $X_{i,k}$ the center of $Q_{i,k}$. We discretize the displacement measure at time $ t_k$ by defining
 for $k\in \N$ \[ \nu_{t_k} = \sum_{i=0}^{N_k^d} \nu^{t_k}(Q_{i,k}) \delta_{X_{i,k}}. \] To complete the construction, consider an optimal transport plan $\pi^k$ for $W(\nu_{t_k},\nu_{t_{k+1}})$ and set
\[ \int_{\R^d} \psi \,d \nu_t = \int_{\R^d} \psi\left(\frac{t- t_{k+1}}{t_k - t_{k+1}} x + \frac{t_k - t}{t_k - t_{k+1}} y\right) \,d \pi^k \quad  \textup{ for } t \in [t_{k+1}, t_k].\] Since $\nu_{t_k}$ and $\nu_{t_{k+1}}$ are discrete so is $\nu_t$ for any $t\in [t_{k+1}, t_k]$. 
Thus, if we set $M(t) = \mathcal{H}^0(\supp(\nu_t))$, there are $t \mapsto \psi_i(t)$ and $t \mapsto X_i(t)$ such that 
\[ \nu_t = \sum_{i=1}^{M(t)} \psi_i(t) \delta_{X_i(t)}. \]
As remarked in \cite{dephilippis2023energy}, the optimal transport problem between two sums of Dirac deltas is a linear programming problem. This means that the transport plan $\pi^k$ can be chosen to be \emph{sparse}, in the sense that we have for $t\in  [t_{k+1}, t_k]$
\begin{equation}
\label{est:active-branches}
    M(t) \leq M(t_k) + M(t_{k+1}).
\end{equation}
By construction, $M(t_k)=\mathcal{H}^0(\supp(\nu_{t_k}))=2^{dk}$, thus for $t\in  [t_{k+1}, t_k]$
\begin{equation}
\label{e:deltasestimate}
    M(t)\lesssim 2^{dk}.
\end{equation}
Finally let us show that $t\mapsto \nu_t$ satisfies $\nu_t \rightharpoonup \mu^-$ as $t \rightarrow 0$.
If $t\in (t_{k+1},t_k)$, then 
\begin{equation*}
    W(\nu_{t}, \mu^-)\leq W(\nu_{t},\nu_{t_{k}})+W(\nu_{t_{k}}, \nu^{t_k})+W(\nu^{t^{k}}, \mu^-).
\end{equation*}
By \eqref{e:geodesic} and  triangle inequality again
\begin{equation*}
    W(\nu_{t}, \nu_{t_k})\leq W(\nu_{t_{k+1}}, \nu_{t_k})\leq W(\nu_{t_{k+1}}, \nu^{t_{k+1}})+W(\nu^{t_{k+1}}, \nu^{t_{k}})+W(\nu^{t_{k}}, \nu_{t_k}).
\end{equation*}
Moreover, since we can construct a competitor for $W(\nu^{t_k}, \nu_{t_k})$ by simply sending all the mass of $\nu^{t_k}$ inside $Q_{i,k}$ to $X_{i,k}$, we have
\begin{equation}
\label{e:wasssquare}
 W(\nu^{t_k}, \nu_{t_k}) \lesssim l_k \sqrt{\Phi}  =r 2^{-k}\sqrt{\Phi}.
\end{equation}
Using  \eqref{e:geodesic} again we thus obtain
\begin{multline*}
    W(\nu_{t}, \mu^-)\lesssim W(\nu_{t_{k+1}}, \nu^{t_{k+1}})+W(\nu^{t_{k+1}}, \nu^{t_{k}})+W(\nu^{t_{k}}, \nu_{t_k}) + W(\nu^{t^{k}}, \mu^-)\lesssim \\ r 2^{-k}\sqrt{\Phi}+\left(\frac{t_{k}-t_{k+1}}{\eps}+\frac{t_k}{\eps}\right)W(\mu^-, \mu^{+})
\end{multline*}
and the conclusion follows letting $k\rightarrow \infty$.

\quad \textbf{\emph{Step 2: Energy estimate for $\nu$.}} Let us estimate the energy $I(\nu)$.
Let us start with the perimeter. By \eqref{e:deltasestimate} and H\"older's inequality we get for $t\in  [t_{k+1}, t_k]$
\begin{equation*}
    \dot{P}(\nu,t)=  \sum_{i=1}^{M(t)} \varphi^{\frac{d-1}{d}}_i(t)  \leq \left(\sum_{i=1}^{M(t)} \varphi_i(t)\right)^\frac{d-1}{d}M(t)^\frac{1}{d}\lesssim \Phi^\frac{d-1}{d}2^k.
\end{equation*}
This implies
\begin{equation}
     \label{perimeter_bound_single_branch-oned}
     \begin{split}
          & P(\nu, (t_{k+1}, t_k))  = \int_{t_{k+1}}^{t_k} \left( \sum_{i=1}^{M(t)} \varphi^{\frac{d-1}{d}}_i(t)  \right)\,dt \lesssim  \Phi^\frac{d-1}{d}2^k (t_{k} - t_{k+1}).  
     \end{split}
\end{equation} 
Summing up all the terms and recalling \eqref{discrete_time_steps_construction} we can bound the perimeter by
\begin{equation}
    \begin{split}
        P(\nu, (0,\eps)) \lesssim \Phi^\frac{d-1}{d} \sum_{k\in \N} (t_{k}-t_{k+1}) 2^k \lesssim \Phi^\frac{d-1}{d}  \eps \sum_{k\in \N} \left(2\delta\right)^k\lesssim  \Phi^\frac{d-1}{d} \eps,
    \end{split}\label{perimeter_constr}
\end{equation}
where in the last passage we used that $\delta<\frac{1}{2}$.
Let us now focus on the kinetic energy. By the definition of $\nu_t$, triangle and Young's inequality we get
\begin{multline*}
E_{kin}(\nu, (t_{k+1}, t_k)) = \left(\frac{1}{t_k- t_{k+1}} \right) W^2(\nu_{t_{k+1}}, \nu_{t_k}) \\
      \leq \frac{1}{t_k - t_{k+1}} \left( W(\nu_{t_{k+1}}, \nu^{t_{k+1}}) + W (\nu^{t_{k+1}}, \nu^{t_{k}})+ W (\nu_{t_{k}}, \nu^{t_{k}})\right)^2 \\
         \lesssim \frac{1}{t_k - t_{k+1}} \left( W^2 (\nu^{t_{k+1}}, \nu^{t_{k}}) +  W^2(\nu_{t_{k+1}}, \nu^{t_{k+1}})+ W^2 (\nu_{t_{k}}, \nu^{t_{k}})\right).
\end{multline*}
 By \eqref{e:geodesic} 
\begin{equation*}
    W^2 (\nu^{t_{k+1}}, \nu^{t_{k}}) =\left(\frac{t_k- t_{k+1}}{\eps}\right)^2 W^2(\mu^-, \mu^+)\leq \left(\frac{t_k- t_{k+1}}{\eps}\right)^2r^2\Phi.
\end{equation*}
Combining this with \eqref{e:wasssquare} yields 
\begin{equation*}
    E_{kin}(\nu, (t_{k+1}, t_k)) \lesssim r^2\Phi\left(\frac{t_k- t_{k+1}}{\eps^2}+ \frac{4^{-k}}{t_k- t_{k+1}} \right).
\end{equation*}
This gives the total estimate 
    \begin{multline}
    \label{e:cineticestimate}
         E_{kin}(\nu , (0,\eps)) \lesssim r^2\Phi \left( \sum_{k\in \N} \frac{t_k - t_{k+1}}{\eps^2}+ \left( \sum_{k\in \N} \frac{1}{4^k (t_k - t_{k+1})}\right) \right)\\
        \lesssim \frac{r^2}{\eps}\Phi \left( \sum_{k\in \N} \delta^k+ \left( \sum_{k\in \N} (4\delta)^{-k}\right) \right)\lesssim \frac{r^2}{\eps}\Phi,
     \end{multline}
where in the last passage we use that $\delta\in(1/4,1/2)$.
Combining \eqref{e:cineticestimate} with \eqref{perimeter_constr} we finally obtain the energy estimate
\begin{equation}\label{e:enerInueps}
    I(\nu, (0, \eps)) \lesssim \frac{r^2}{\eps}\Phi +  \Phi^\frac{d-1}{d} \eps.
\end{equation}
\quad \textbf{\emph{Step 3: General construction.}}
Let us now construct a competitor for generic $x_\pm$. For $x_0\in \R^d$, let $\tau_{x_0}(x)= x-x_0$. Set $h=x_+-x_-$ and consider
\begin{equation*}
    \nu^-=\tau_{h}\#\mu^-.
\end{equation*}
We can apply Step $1$ to $\nu^-,\mu^+$ and obtain a competitor $\nu\in \A^{*}(\R^d,(0,\eps))$ such that $\nu_0=\nu^-$, $\nu_\eps=\mu^+$ and \eqref{e:localconstruction} holds. 
For $t\in(0,\eps)$, we then set
\begin{equation*}
    \mu_t=\tau_{-\frac{\eps-t}{\eps}h}\#\nu_t.
\end{equation*}
Notice that $\mu\in \A^{*}(\R^d,(0,\eps))$,
\begin{equation*}
   \mu_0=\tau_{-h}\#\nu_0=\mu^- \quad \text{ and } \quad  \mu_\eps=\tau_{0}\#\nu_\eps=\mu^+.
\end{equation*}
Moreover (see Figure \ref{fig:construction_single_branch}),
\begin{equation*}
    \supp(\mu_t) \subseteq \frac{\eps - t}{\eps} B_{r}(x_-) + \frac{t}{\eps} B_{r}(x_+) \quad \textup{ for all } t \in [0, \eps].
\end{equation*}
\begin{figure}[h]
    \centering\includegraphics{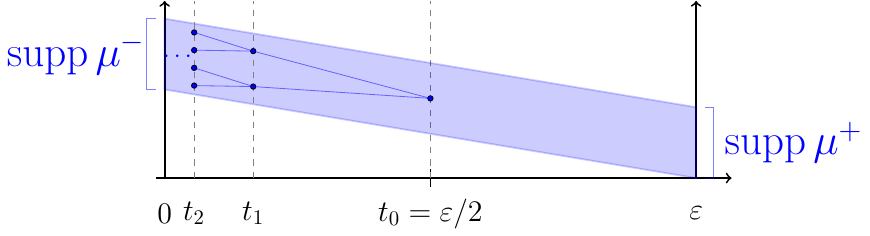}
    \caption{Illustration of the construction on one single branch for $d=1$. The support of $\mu_t$ is contained in the blue domain.}
    \label{fig:construction_single_branch}
\end{figure}

Let us estimate the energy of $\mu$. Clearly 
\begin{equation}
\label{e:pertranslate}
    P(\mu,(0,\eps))=P(\nu,(0,\eps)).
\end{equation}
Let us now compute the kinetic energy. If $\mu_t = \sum_{i} \varphi_i(t) \delta_{\Tilde{X}_i(t)}$ and $\nu_t = \sum_{i} \varphi_i(t) \delta_{X_i(t)}$, then by construction
\begin{equation*}
    \Tilde{X}_i(t)=X_i(t)-\frac{\eps-t}{\eps}h.
\end{equation*}
This implies 
\begin{equation*}
    \abs{\dot{\Tilde{X}}_i(t)}^2=\abs{\dot{X}_{i}(t)}^2+\frac{\abs{h}^2}{\eps^2}+\frac{2}{\eps}h\cdot \dot{X}_i(t),
\end{equation*}
which gives 
\begin{equation}
\label{e:kintranslate}
    E_{kin}(\mu , (0,\eps))=E_{kin}(\nu , (0,\eps))+\frac{\abs{h}^2}{\eps}\Phi+\frac{2}{\eps}h\cdot \int_{0}^\eps\sum_i\varphi_i\dot{X}_i(t)\,dt.
\end{equation}
Take $\pi$ to be an optimal transport plan for $W(\mu^-,\mu^+)$ and let $F(x,y)=(x-h,y)$. Then, $\pi$ can be chosen such that 
\begin{equation*}
    \pi=F\#\overline{\pi}
\end{equation*}
where $\overline{\pi}$ is the optimal transport plan for $W(\nu^-,\mu^+)$.
This implies that 

\begin{multline}
\label{e:splitwass}
    W^2 (\mu^-,\mu^+)= \int_{\R^d\times \R^d}\abs{x-y}^2\,d\pi = \int_{\R^d\times \R^d}\abs{x-h-y}^2\,d\overline{\pi}\\=  W^2 (\nu^-,\mu^+)+ \abs{h}^2\Phi+2h\cdot \int_{\R^d\times \R^d}(y-x)\,d\overline{\pi}.
\end{multline}
We claim that 
\begin{equation}
\label{e:equivalencetransport}
     \int_{0}^\eps\sum_i\varphi_i\dot{X}_i\,dt=\int_{\R^d\times \R^d}(y-x)\,d\overline{\pi}.
\end{equation}
Then by \eqref{e:pertranslate},\eqref{e:kintranslate},\eqref{e:splitwass} and \eqref{e:equivalencetransport}, we would get 
\begin{equation*}
    I(\mu)= I(\nu)+\frac{1}{\eps}(W^2 (\mu^-,\mu^+) -W^2 (\nu^-,\mu^+))\leq I(\nu)+\frac{1}{\eps} W^2 (\mu^-,\mu^+).
\end{equation*}
Combining this with \eqref{e:enerInueps} would conclude the proof of \eqref{e:localconstruction}.
We thus prove \eqref{e:equivalencetransport}. The computation is similar to \cite[Lemma 3.1]{minimizers2d}. By definition,
\begin{equation*}
    \int_{\R^d\times \R^d}(y-x)\,d\overline{\pi}=\int_{\R^d} y\,d\mu^+-\int_{\R^d} x\,d\nu^-.
\end{equation*}
Let now $\delta<\frac{\eps}{2}$ and split
\begin{equation}
\label{e:splitcompute}
    \int_{0}^\eps\sum_i\varphi_i\dot{X}_i\,dt=\int_{0}^\delta\sum_i\varphi_i\dot{X}_i\,dt+\int_{\delta}^{\eps-\delta}\sum_i\varphi_i\dot{X}_i\,dt+\int_{\eps-\delta}^\eps\sum_i\varphi_i\dot{X}_i\,dt.
\end{equation}
In $[\delta,\eps-\delta]$, $\nu$ has by construction a finite number of branches, hence applying \cite[Lemma 5.9]{CGOS}, see also \cite[Lemma 2.6]{minimizers2d}, we find
\begin{equation*}
    \int_{\delta}^{\eps-\delta}\sum_i\varphi_i\dot{X}_i(t)\,dt=\sum_i\varphi_i(\eps-\delta)X_i(\eps-\delta)-\sum_i\varphi_i(\delta)X_i(\delta)=\int_{\R^d}y\,d\mu_{\eps-\delta}-\int_{\R^d}x\,d\mu_\delta.
\end{equation*}
Sending $\delta\to 0$, we find
\begin{equation*}
   \lim_{\delta\to 0} \int_{\delta}^{\eps-\delta}\sum_i\varphi_i\dot{X}_i(t)\,dt= \int_{\R^d}y\,d\mu^+-\int_{\R^d}x\,d\nu^-.
\end{equation*}
To conclude, let us show that the other terms in \eqref{e:splitcompute} are vanishing when $\delta\to 0$. By H\"older inequality applied  twice and since $\sum_i\varphi_i=1$
\begin{equation*}
    \int_{0}^\delta\sum_i\varphi_i\dot{X}_i\,dt\leq \int_{0}^\delta\left(\sum_i\varphi_i\abs{\dot{X}_i}^2\right)^\frac{1}{2}\,dt\leq \delta^\frac{1}{2}\left(\int_{0}^\delta\sum_i\varphi_i|\dot{X}_i|^2\,dt\right)^\frac{1}{2} \stackrel{\delta\to 0}{\rightarrow} 0.
\end{equation*}
 Arguing similarly on the interval $[\eps-\delta,\eps]$ concludes the proof.
\end{proof}
\subsection{Local scaling laws for concentrated boundary measures $\mu_T$}

In this section we fix $d=1$.

\begin{proof}[Proof of Theorem \ref{thm:local_energy_lower}]
Let $r\in (0,1)$. We start by covering the support of $\mu_T$ with pairwise disjoint intervals of radius proportional to $r$ . Using the uniform lower bound \eqref{e:malphacon}, we have that the number of these intervals is at most of the order of $r^{-\alpha}$. Next, by the cone property, at $t = T- \eps$  the set of points that irrigate one of these intervals can be decomposed into at most three intervals of length proportional to $r$. Having the decomposition at $t=T-\eps$ and $t=T$ we construct the competitor with respect to each of the intervals on the left boundary by applying Lemma \ref{lem:single_branch}.\\

\quad \emph{\textbf{Step 1. Decomposition of left and right boundaries.} } Since $\mu_T$ is lower $\alpha-$Ahlfors regular,  for every $x\in \supp(\mu_T)$ and $r\in (0,1)$ we have  $\mu_T(B_r(x))\gtrsim_M r^\alpha$. Let us fix $r\in (0,1)$. By Vitali's covering theorem \cite[Theorem 1.24]{evans2018measure}, there exists $x_i\in \supp(\mu_T)$ for $i=1,\dots, N$ such that $B_r(x_i)\cap B_r(x_j)= \emptyset$ for $i\neq j$ and $\supp(\mu_T) \subset \cup_{i=1}^N B_{5r}(x_i) $.
By  hypothesis we have
\begin{equation}\label{e:numberballs}
    1=\mu_T(\Tor)\geq \mu_T(\cup_{i=1}^N B_{r}(x_i)) = \sum_{i=1}^N \mu_T(B_r(x_i))\gtrsim_M  Nr^\alpha.
\end{equation}
 Since we are in dimension $1$, we can increase the radius of the balls $B_r(x_i)$ until they either touch or reach the radius $5r$. We thus find a collection of at most $N$ disjoint intervals $J_i(T)$ of length at most $10r$ that cover $\supp(\mu_T)$. Consider the parametric representation $X(x,t)$ of $\mu$ given by Proposition \ref{prop:lagrangian}. By Corollary \ref{cor:3intervals}, for each of the intervals $J_i$ and for all $t\in[T-\eps,T]$ there exist at most three intervals $J_i^j(t)$ such that for $\mu_T$-a.e. $ x \in J_i, \, X(x,t)\in J_i^j(t)$ for some $j=1,2,3$, see Figure \ref{fig:dec_single}. Moreover, the three intervals satisfy $|J_i^j(t)|\leq |J_i|\leq 10r$ for $j=1,2,3$ .  \\

\begin{figure}[h]
    \centering\includegraphics[scale = 0.5]{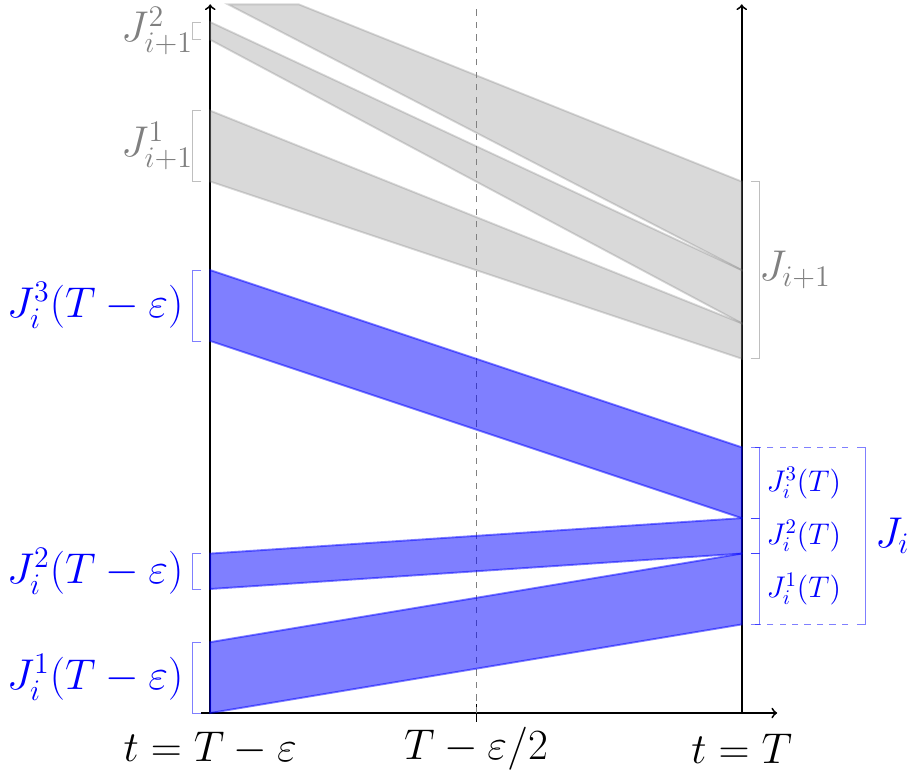}
    \caption{Sketch of the described decomposition with respect to $J_i$. The part of the minimizer which irrigates $J_i$ is contained in the blue domain.}
    \label{fig:dec_single}
\end{figure} 

\quad \emph{\textbf{Step 2. Construction on $J_i^j$}}

We introduce the restricted measure $\mu^{i,j}_T$ by setting
\[ \mu^{i,j}_T  = \mu_T\restr J_i^j(T) \quad \textup{ for }    i = 1, \dots, N \textup{ and } j=1,2,3. \]
Let us define
\begin{equation*}
     \mu^{i,j}_{T-\eps}  = X(\cdot,T-\eps)\#\mu^{i,j}_T \quad \textup{for }   i = 1, \dots, N \textup{ and } j=1,2,3.
\end{equation*}
In the language of Proposition \ref{prop:subsys}, this is the backward subsystem starting from $\mu^{i,j}_T$. Since the supports of $\mu^{i,j}_T$ are disjoint,
\begin{equation*}
    \sum_{i,j} \mu^{i,j}_{T-\eps}=\mu_{T-\eps}.
\end{equation*}
Moreover, by the decomposition in Step 1, we have 
\begin{equation} 
    \label{supp:optimal_plan}
    \supp( \mu_{T-\eps}^{i,j}) \subseteq  J_i^j(T-\eps).
\end{equation}
Hence, for $t\in \{T-\eps, T\}$
\[|\supp(\mu_t^{i,j})| \leq |J_i^j(t)| \lesssim r.\]
Finally, notice that since by Lemma \ref{lem:monotonecoupling}, $X(\cdot, T-\eps)$ is monotone (see also Remark \ref{rem:periodicnonperiodic}), it is an optimal transport map for both $W^2(\mu_T,\mu_{T-\eps})$ and $W^2( \mu_T^{i,j},\mu_{T-\eps}^{i,j})$ so that
\begin{equation}
\label{e:wassersteinsplit}
    \sum_{i,j} W^2(\mu_{T-\eps}^{i,j}, \mu_T^{i,j})= W^2(\mu_{T-\eps}, \mu_T).
\end{equation}
 We now apply the previous lemma within the cones, i.e. the blue marked domains in Figure \ref{fig:construction}.
\begin{figure}[h]
    \centering
    \includegraphics[scale = 0.7]{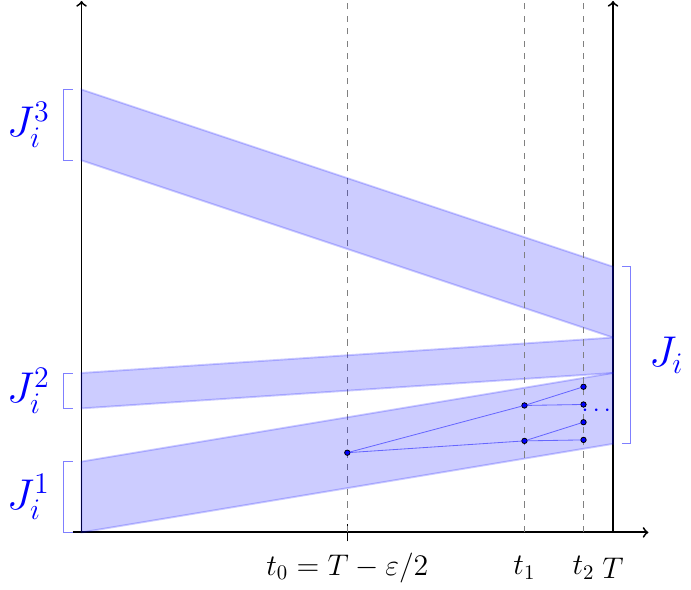}
    \caption{Sketch of the local construction within the cones.}
    \label{fig:construction}
\end{figure}

By Lemma \ref{lem:single_branch}, for every $i,j$, there is an admissible construction $\nu^{i, j}\in \mathcal{A}^*(\Tor, (T-\eps, T))$ such that $\nu^{i, j}_T = \mu^{i, j}_T$, $\nu^{i, j}_{T-\eps} = \mu^{i, j}_{T-\eps}$ and
\begin{equation}
    \label{est:energy_Ji}
    I(\nu^{i,j}, (T-\eps,T)) \leq \frac{1}{\eps} W^2(\mu_{T-\eps}^{i,j}, \mu_T^{i,j}) + C\left( \frac{r^2}{\eps}\mu_T(J_i^j(T)) + \eps\right).
\end{equation}
Notice that by \eqref{e:supmut}, the supports of $\nu^{i,j}$ are pairwise disjoint. We finally set $\nu =  \sum_{i,j} \nu^{i,j}.$ We have  $\nu \in \mathcal{A}^*(\Tor,(T-\eps, T))$,  $\nu_T = \mu_T$ and $\nu_{T-\eps} = \mu_{T-\eps}$. The minimality of $\mu$ together with estimates \eqref{est:energy_Ji}, \eqref{e:wassersteinsplit} and subadditivity of the energy imply that for some $C=C(M)>0$,
    \begin{multline*}
        I(\mu, (T-\eps,T)\leq I(\nu, (T-\eps,T)) \leq \sum_{i,j} I(\nu^i, (T-\eps, T))\\ \leq  \frac{1}{\eps} W^2(\mu_{T-\eps}, \mu_{T}) + C\left(  \frac{ r^2}{ \eps } \sum_{i,j}\mu_T(J_i^j(T)) + N \eps \right)
         \stackrel{\eqref{e:numberballs}}{\leq} \frac{1}{\eps} W^2(\mu_{T-\eps}, \mu_T) +  C\left( \frac{ r^2}{ \eps } + r^{-\alpha}\eps\right).
      \end{multline*}
 \\
\quad \emph{\textbf{Step 3. Equipartition of energy and conclusion} .}
 By the previous step,
    \begin{equation*}\label{est:local_energy_final}
    \begin{split}
        I(\mu, (T-\eps,T))\leq \frac{1}{\eps} W^2(\mu_{T-\eps}, \mu_T) +  C \left(\frac{ r^2}{ \eps } + r^{-\alpha}\eps \right).
    \end{split}
\end{equation*}
Choosing $r=  \eps^{2/(2+\alpha)}\in(0,1) $  to optimize the right-hand side, we  obtain
    \begin{equation}\label{est:local_energy_final_optimized}
    \begin{split}
        I(\mu, (T-\eps,T))\leq \frac{1}{\eps} W^2(\mu_{T-\eps}, \mu_T) + C \eps^{\frac{2-\alpha}{2+\alpha}}.
    \end{split}
\end{equation}
Since \eqref{e:ecinestimate} holds, we obtain 
\begin{equation*}
   P(\mu, (T-\eps, T)) \lesssim_M  \eps^{\frac{2-\alpha}{2+\alpha}}.
\end{equation*}
Finally, by the equipartition of energy, see Proposition \ref{prop:equipartition}, there is a constant $\Lambda$ with $\Lambda \leq I(\mu, (0, \eps))/T$ such that
\begin{equation}
    P(\mu, (T-\eps, T)) = E_{kin}(\mu, (T-\eps, T)) +\Lambda \eps \label{identity:perimeter},
\end{equation}
from which the claimed \eqref{e:upperboundth1.6} follows.
\end{proof}

\section*{Acknowledgement}
The authors gratefully acknowledge funding of the Deutsche Forschungsgemeinschaft (DFG, German
Research Foundation) via project 195170736 - TRR 109 and the FMJH program Junior Scientific Visibility. A.C. is partially supported by the European Union's Horizon 2020 research and innovation programme under the Marie Sklodowska-Curie grant agreement No 945332. The authors would like to thank Martin Huesmann for his useful remarks on Proposition \ref{lem:reg_displacement}.
\flushright \includegraphics[height=.9cm]{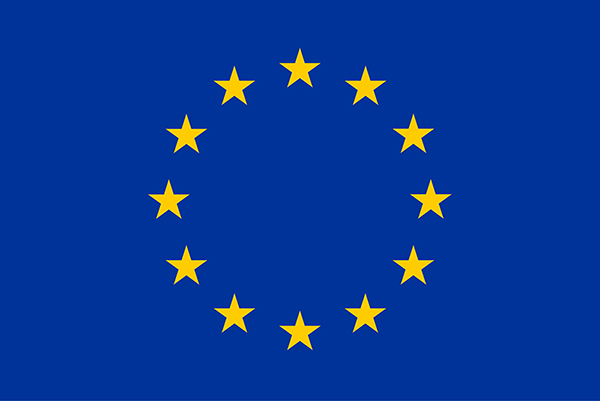}

\printbibliography

\medskip
\small
\begin{flushright}
\noindent \verb"acosenza@math.univ-paris-diderot.fr"\\
Universit\'e Paris Cit\'e and Sorbonne Universit\'e,\\
CNRS, Laboratoire Jacques-Louis Lions (LJLL),\\
F-75006 Paris, France
\end{flushright}
\medskip
\small
\begin{flushright}
\noindent \verb"michael.goldman@cnrs.fr"\\
CMAP, CNRS, \'Ecole polytechnique, Institut Polytechnique de Paris,\\ 
91120 Palaiseau, France
\end{flushright}
\medskip
\begin{flushright}
\noindent \verb"melanie.koser@hu-berlin.de"\\
Humboldt-Universit\"at zu Berlin, Institut f\"ur Mathematik,\\ 
10099, Berlin, Germany
\end{flushright}

\end{document}